\documentclass[12pt]{article}

\begin{document}

\input epsf.sty

\centerline{\bf S.P.Novikov\footnote{Sergey P. Novikov, IPST and
MATH Department, University of Maryland, College Park MD, USA and
Landau Institute, Moscow, e-mail novikov@ipst.umd.edu; This work
is partially supported by the Russian Grant in the Nonlinear
Dynamics }}

 \vspace{0.3cm}

\centerline{\Large  Topology of Foliations}
 \centerline{\Large given by the real part of holomorphic 1-forms}

\vspace{0.5cm}

{\it Abstract.
  Topology of Foliations of the Riemann Surfaces given by the
 real part of
generic  holomorphic 1-forms, is studied. Our approach is based on
the notion of Transversal Canonical Basis of  Cycles (TCB) instead
of using just  one closed transversal curve as in the classical
approach of the ergodic theory. In some cases the TCB approach
allows us to present a convenient combinatorial model of the whole
topology of the flow, especially effective for g=2. A maximal
abelian covering over the Riemann Surface provided by the Abel
Map, plays a key role in this work. The behavior of our system in
the Fundamental Domain of that covering can be easily described in
the sphere with $g$ holes. It leads to the Plane Diagram of our
system.  The complete combinatorial model of the flow is
constructed. It is based on the Plane Diagram and g straight line
flows in the planes corresponding to the g canonically adjoint
pairs of cycles in the Transversal Canonical Basis. These pairs do
not cross each other. Making cuts along them, we  come to the
maximal abelian fundamental domain (associated with Abel Map and
Theta-functions) instead of the standard 4g-gon in the Hyperbolic
Plane and its beautiful ''flat'' analogs which people used for the
study of  geodesics of the flat metrics with singularities. }

\vspace{0.5cm}

{\bf Introduction}. The family of parallel straight lines in the
Euclidean Plane $R^2$ gives us after factorization by the lattice
$Z^2\subset R^2$ the standard straight line flow in the 2-torus
$T^2$. It is a simplest ergodic system for the  irrational
direction. This system is Hamiltonian with multivalued Hamiltonian
function $H$ and standard canonically adjoint euclidean
coordinates $x,y$ (i.e. the  1-form $dH$ is closed but not exact,
and $x_t=H_y, y_t=-H_x$ where $H_x,H_y$ are constant). Every
smooth Hamiltonian system on the 2-torus without critical points
with irrational ''rotation number'' is diffeomorphic to the
straight line flow. Every $C^2$-smooth dynamical system on the
2-torus without critical points and with irrational rotation
number is $C^0$-homeomorphic to the straight line flow according
to the famous classical theorem (it is not true for $C^1$).

{\bf Question: What is a right analog of the straight line flow
for the Riemann Surfaces of higher genus?}

 Certainly, it should
be some special class of Hamiltonian Systems with multivalued
Hamiltonian and trajectories given by foliation $dH=0$. Which
special subclass has the best properties? For the genus $g=2$ the
answer to this question will be given in the last section of this
article.

 In
many cases we ignore time dependence of the trajectories and
discuss only the properties of  foliation $dH=0$ given by the
closed 1-form. The present author started to investigate such
foliations  in early 1980s as a part of the newborn Topology of
The Closed 1-Forms. An important example was found in the Quantum
Solid State Physics describing the  motion of  semiclassical
electrons along the so-called Fermi Surface for the single crystal
normal metals and low temperature in the strong magnetic field
(see \cite{N82}). Good introduction to the corresponding Physics
can be found  in the textbook \cite{Abr}. The Fundamental
''Geometric Strong Magnetic Field Priciple'' was formulated by the
Kharkov school of I.Lifshitz many years ago (and fully accepted by
physics community). It claims that following geometric picture
gives a good description of  electrical conductivity in the single
crystal normal metals in the ''reasonably strong'' magnetic field
(of the size $1t<B<10^3t$ and temperature $T<1K$ for the normal
metal like gold): There is a Fermi Surface $M_F$ in the 3-torus of
quantum ''quasimomenta'' $p\in T^3$ where $M_F\subset T^3$ is a
nondegenerate level $\epsilon=\epsilon_F$ of the Morse function
$\epsilon :T^3\rightarrow R$ called the ''dispersion relation''.
Every constant homogeneous  magnetic field $B$ defines a 1-form
$\sum_i B_idp_i$ whose restriction to the Fermi Surface is exactly
our closed 1-form: $dH=\sum_i B_idp_i|_{M_F}$. The electron
trajectories exactly coincide with connectivity components of the
sections of Fermi Surface by the planes orthogonal to magnetic
field. One might say that they are the levels of quasiperiodic
function on the plane with 3 periods. An extensive study of that
class was performed by the author's Moscow Topological Seminar
since early 1980s. These studies were continued in Maryland since
the second half of 1990s. Among the author's students who made
important contribution here, let me mention A.Zorich, S.Tsarev,
I.Dynnikov, A.Maltsev, R.Deleo (see the survey articles \cite{D99,
MN04,DN05} describing
 topological, dynamical and physical results of our studies).
This class of systems has remarkable ''topological complete
integrability''  and ''topological resonance'' properties in the
nonstandard sense, for the set of directions of magnetic field of
the full measure on the 2-sphere $S^2$. These properties play a
key role in the physical applications.  Numerical studies are
described in \cite{RdL}. According to the fundamental results of
that theory, topological properties here generically can be
reduced  to the case of the genus 1. (Let us mention that the
hamiltonian systems on  2-torus were studied in \cite{A}; their
ergodic properties were found finally in \cite{KhS}). Our problem
is more complicated. In particular, there exists a nonempty set of
parameters with Hausdorf Dimension $d\leq 1$ (presumably, even
strictly less than one) on the 2-sphere, where a complicated
''chaotic'' behavior was discovered.

 Following  class of systems was proposed  by the
various mathematicians (no applications outside of pure
mathematics were found for them  until now, unfortunately): Take
any nonsingular compact Riemann Surface $V$, i.e complex algebraic
curve, with genus equal to $g$. There exists a $g$-dimensional
complex linear space $C^g$ of holomorphic differential forms with
basis $\omega_1,\dots,\omega_g\in C^g$. Every holomorphic 1-form
$\omega\in C^g$ defines a Hamiltonian system (foliation) $\Re$ on
the manifold $V$: $$\Re=\{ \omega^R=0\}$$ because $d(\omega^R)=0$.
Here $\omega=\omega^R+i\omega^I$. There is even more general class
of ''foliations $\Re$ with invariant transversal measure'' on
Riemann Surfaces: the existence  of  measure is required on the
intervals transversal to the leaves (trajectories) invariant under
the deformations such that every point is moving along
trajectories. Take any holomorphic quadratic differential $\Omega$
and define foliation $\Re$ by the formula $(\sqrt{\Omega})^R=0$.
This is a locally hamiltonian foliation $\Re$ (i.e. it admits a
transversal measure) but may be non-orientable (it does not admit
 time direction globally), so we do not consider them.

The  systems with transversal invariant measure were studied since
early 1960s by the following method (see in the book \cite{KHa}):
Take any closed curve $\gamma$ transversal to our foliation $\Re$.
Assume that almost every nonsingular trajectory is dense. For
every point $Q\in \gamma$ except finite number there exists a
first time $t_Q$ such that trajectory started in the point $Q$
returns to the curve $\gamma$ as a new point $P$ (in positive
direction of time). We define a ''Poincare map'' $Q\rightarrow P$.
This map by definition preserves a transversal measure which is a
restriction of the form $dH$ on the curve $\gamma$.   {\bf It is
the Energy Conservation Law} for  Hamiltonian system. So our
transversal closed curve $\gamma$ is divided into
 $k=k_{\gamma}$ intervals $\gamma=I_1+I_2+\ldots+I_k$. The
Poincare map looks here as a simple permutation of these intervals
on the circle; it is ill-defined in the finite number of points
only. The time $t_Q$ varies continuously within each interval
$Q\in I_j$.

No doubt, the use of  closed transversal curves is extremely
productive. At the same time,  we are not satisfied by this
approach; Following questions can be naturally asked:

1.This method essentially ignores time and length of trajectories
starting and ending in $\gamma$. Where they are traveling and how
long? We would like to see  some sort of global topological
description  of the flow (or foliation $\Re$) on the algebraic
curve $V$ {\bf and its abelian coverings}  similar to the case of
genus 1 as much as possible. The ergodic  characteristics of
foliation  can be found here, but this model certainly in not
enough for the description of topology.

 2.There are many different closed transversal curves
in the foliation $\Re$. How this picture depends on the choice of
tranversal curve $\gamma$? Nobody classified them yet as far as I
know. Indeed, in the theory of codimension 1 foliations developed
by the present author in 1960s (see\cite{N65}) several algebraic
structures were defined for the closed transversal curves:

Fix any nonsingular point $Q\in V$ and consider all closed
positively (negatively) oriented transversal curves starting and
ending in $Q$. We can multiply them. Transversal homotopy classes
of such curves generate {\bf A Transversal Semigroup}
$\pi_1^+(\Re,Q)$ and its natural homomorphism into the fundamental
group (even, into the fundamental group of the unit tangent
$S^1$-bundle $L(V)$
$$\psi^{\pm}:\pi_1^{\pm}(\Re,Q)\rightarrow\pi_1(L(V),Q)\rightarrow\pi_1(V,Q)$$
The set of all closed transversal curves naturally maps into the
set of conjugacy classes in fundamental group. We denote it also
by $\psi$. The Transversal Semigroups might depend of the leaf
where the initial point $Q$ is chosen. For example, they are
different for the separatrices entering critical points and for
the generic nonsingular leaves.

{\bf How to calculate these invariants for the Hamiltonian
foliations described above on the algebraic curves?}

Let me point out the simplest fundamental properties of these
foliations:

Property 1. They have only saddle type critical points. In the
generic case such foliation has exactly 2g-2 nondegenerate
saddles.

Property 2. Every nonempty closed  transversal curve $\gamma$ is
non-homologous to zero. The period of the form $dH$ is positive
$\oint_{\gamma}dH>0$ for every positively oriented transversal
curve. So the composition $$\pi_1^+(\Re,Q)\rightarrow
\pi_1(V,Q)\rightarrow H_1(V,Z)\rightarrow R$$ does not map any
element into zero. Its image is strictly positive.

 We present here a theory of the special classes {\bf $T$} and {\bf
$T^{k}$} of these foliations. Some of them can be considered (in
the generic case, at least) as a natural higher genus analog of
the straight line flow, especially for $g=2$. The definition of
the classes is following:

{\bf Definitions}

a.We say that foliation on the Riemann surface $V$   belongs to
the {\bf Class} $T$  if there exist a {\bf Transversal Canonical
Basis} of curves $$a_1,b_1,\ldots,a_g,b_g$$ such that

I.All these curves are nonselfintersecting and {\bf transversal}
to the foliation.

II.The curves $a_j$ and $b_j$ cross each other transversally for
all $j=1,\ldots,g$ exactly in one point. Other curves do not cross
each other.

b.We say that foliation  belongs to the {\bf Class} $T^{0}$  if
there exists a canonical basis such that all $a$-cycles
$a_1,\ldots,a_g$ are non-selfintersecting, do not cross each other
and are transversal to the foliation .

c.We say that foliation  belongs to the {\bf Mixed Class}  of the
Type $T^{k}$  if there exists {\bf an incomplete canonical basis}
$a_j,b_q,j\leq g,q\leq k,$ such that all these cycles are
non-selfintersecting. The only pairs crossing each other are $a_j$
and $b_j$ for $j=1,2,\ldots,k$. The intersections of $a_j$ and
$b_j$ are transversal and consist of one point each. All these
curves $a_j,b_q,j\leq g,q\leq k$ are transversal to foliation.  We
have $T=T^{g}$

\newtheorem{rem}{Remark}

\begin{rem}
A number of people including Katok, Hasselblatt, Hubbard, Mazur,
Veech,
 Zorich, Konzevich, McMullen  and others wrote a lot of works related to study the
ergodic properties  of foliations with  ''transversal measure'' on
the Riemann Surfaces, and their total moduli space
 (see  \cite{K,Sat,KHa,HM,KMS,M,V,MT,Z,Z1,Z2,ZK,Kn,F,McM}).  People investigated
recently  closed geodesics of the flat Riemannian Metric
$ds^2=\omega\bar{\omega}$ singular in the critical points of the
holomorphic 1-form $\omega$. These geodesics consists of all
trajectories of the Hamiltonian systems
$(\exp\{i\phi\}\omega)^R=0$ on the algebraic curve $V$.  Closed
geodesics appear for the special  number of angles $\phi_n$ only.
Therefore we don't see them describing the generic Hamiltonian
systems of that type. Beautiful analogs of the Poincare' 4-gons
associated with these flat metrics
 with singularities were invented  and used.
 Our intention is to describe completely
topology of some specific good classes of the generic foliations
given by the equation $\omega^R=0$ for the holomorphic form
$\omega$. As A.Zorich pointed out to the author, considering very
specific examples people certainly observed some features which
may illustrate our ideas (see
 \cite{KHa,McM}).
 It seems that our key idea of transversal canonical
 basis did not  appeared before.
\end{rem}

\vspace{0.3cm}

 \centerline{\Large Section 1. }

\centerline{\Large Definitions. The Case of Hyperelliptic Curves.}

\vspace{0.5cm}

First of all, our goal is to show the examples of concrete
foliations in the classes $T,T^{k}$. Consider  any real
hyperelliptic curve of the form
$$w^2=R_{2g+2}(z)=\prod_j (z-z_j),z_j\neq z_l$$ where all roots
$z_j$ are real and ordered naturally $z_1<z_2<\ldots <z_{2g+2}$.

\newtheorem{lem}{Lemma}
\begin{lem}
For every complex non-real and non-imaginary number $u+iv,u\neq
0,v\neq 0$, every polynomial $P_{g-1}(z)$ and  number $\epsilon$
small enough, the Hamiltonian system defined by the closed
harmonic 1-form below belongs to the class  $T^{0}$ with  possible
choice of $a$-cycles as any subset of $g$ cycles out of
$a_1,\ldots,a_g$ and $c_1,\ldots,c_{g+1}$ not crossing each other:
$$ \omega=\frac{(u+iv+\epsilon P_{g-1}(z))dz}{\sqrt{R_{2g+2}(z)}},
\omega^R=0$$

 For $g=1,2$ this foliation belongs to the class $T$.
For $g>2$ it belongs to the class $T^{2}$.

\end{lem}
Proof. The polynomial $R=R_{2g+2}(z)$ is real on the cycles
$a_j=p^{-1}[z_{2j}z_{2j+1}],j=1,...,g$ and purely imaginary on the
cycles $c_q$ located on the real line $x$ immediately before and
after them $c_q=p^{-1}[z_{2q-1}z_{2q}],q=1,2,\ldots,g+1$.

 For $g=1$ we take a canonical basis
$a_1,b_1=c_1$. For $g=2$ we take a canonical basis
$a_1,a_2,b_1=c_1,b_2=c_3$. We shall see that all these cycles are
transversal to foliation.

  For
$\epsilon=0$ we have $$\omega^R=(udx-vdy)/R^{1/2}=0,y=0,x\in a_j$$
$$\omega^R=(vdx+udy)/iR^{1/2}=0, y=0,x\in c_j$$ In both cases we
have for the second component $dy\neq 0$ for the direction of
Hamiltonian system. So the transversality holds for all points on
the cycles except (maybe) of the branching points $z_j$. The
cycles $c_l,a_k$ are orthogonal to each other in the crossing
points (i.e. in the branching points).  We need to check now that
in all branching points $z_j$ the angle between trajectory and
both cycles $a,c$ is never equal to $\pm\pi/2$. because they meet
each other with this angle exactly. After substitution
$w'^2=z-z_j$ we have
$$\omega=(u+iv)2w'dw'/w'F_j^{1/2}=2(u+iv)dw'/F^{1/2}_j$$ where
$F_j=\prod_{l\neq j} (z-z_l)$. We see that the real function
$F_j(x)$ for real $x\in R$ near the point $z_j\in R$, does dot
change sign passing through this point $x=z_j$:  $F_j(z_j)\neq 0$.
For $w'=f+ig$ we have $\omega^R=2(udf-vdg)/F_j^{1/2}$ if
$F_j(z_j)>0$ or $\omega^R=2i(udg+vdf)/F^{1/2}$ if $F_j(z_j)<0$. We
have $x=f^2-g^2,y=2fg$. The condition $y=0$ implies the equation
$fg=0$, i.e. the union of the equations $f=0$ and $g=0$. It is
exactly an orthogonal crossing of our cycles. In both cases our
system $\omega^R=0$ implies transversality of trajectories to the
both cycles $df=0$ and $dg=0$ locally.

We choose now the $a$-cycles as it was indicated above. All of
them are transversal to foliation. Therefore we are coming to the
class $T^0$. Now we choose following two $b$-cycles: $b_1=c_1$ and
$b_g=c_{g+1}$. According to our arguments, this choice leads to
the statement that our foliation belongs to the class $T$ for
$g=1,2$ and  to the class $T^2$ for all $g\geq 2$.

 Lemma
is proved now because  small $\epsilon$-perturbation cannot
destroy transversality along  the finite family of compact cycles.

We choose now the generic perturbation such that {\bf all critical
points became nondegenerate, and all saddle connections and
periodic trajectories nonhomologous to zero, disappear.}

Let now $R=R_{2g+2}=\prod_{j=1}^{2g+2}(z-z_j)$ is a polynomial of
even degree as above with real simple roots $z_j\in R$, and
$\omega=P_{g-1}(z)dz/R^{1/2}$ is a generic holomorphic 1-form. Let
$P=u+iv$ where $u,v$ are real polynomials in the variables $x,y$.
Only the zeroes $v=0$ in the segments
$[z_{2q-1}z_{2q}],q=1,...,g+1$, and the zeroes $u=0$ in the
segments $[z_{2j}z_{2j+1}],j=1,...,g$, are important now. We
assume that $u(z_k)\neq 0$ and $v(z_k)\neq 0$ for all
$k=1,...,2g+2$.
\begin{lem}  Remove all open segments  containing the important zeroes, i.e. the
open segments $[z_{2q-1}z_{2q}]$ containing the zeroes  $v=0$, and
all open segments $[z_{2j}z_{2j+1}]$ containing the zeroes $u=0$.
If remaining segments are enough for the construction of the
half-basis $a_1,\ldots,a_g$ (i.e. there exists at least $g$
disjoint closed segments between them), then our foliation belongs
to the class $T^0$. In particular, it is always true for $g=2$
where $u$ and $v$ are the linear functions (and have no more than
one real zero each).

 Let there are no ''real'' (i.e. located on the $x$-line) zeroes $v=0$
  in the segments $[z_1z_2],[z_5z_6]$,
and no real zeroes $u=0$ in the segments $[z_2z_3],[z_4z_5]$  for
$g=2$; Then this foliation belongs to the class $T$.

 Let $g>2$, all real   zeroes of polynomial $v=0$ belong to the open intervals
$$(-\infty,z_1),(z_3z_4),(z_5z_6),\ldots,
 (z_{2g-1}z_{2g}),( z_{2g+2},+\infty)$$, and all real zeroes $u=0$
are located  on the $x$-line in the open segments
$$(-\infty,z_2),(z_3z_4),...,(z_{2g-1}z_{2g},(z_{2g+2},+\infty)$$.
In this case  foliation $\omega^R=0$ belongs to the class $T^2$.

\end{lem}
Proof is exactly the same as above.

For every foliation of the class $T^0$ we cut Riemann Surface $V$
along the transversal curves $a_j$. The remaining manifold
$\tilde{V}$ has a boundary
$$\partial{\tilde{V}}=\bigcup_ja_j^{\pm}$$ where $a_j^{\pm}=S^1$
with foliation entering it from inside for $(a^-_j)$ and leaving
it towards the inside of the surface $\tilde{V}$ for $(a_j^+)$.
This system can be considered as a system on the plane in the
domain $\tilde{V}=D^2_*$ where star means that
 $2g-1$ holes removed from this disk inside: The external boundary
is taken as $\partial{D^2_*}=a_1^+$. The boundaries of inner holes
are
 $$a_1^-,a_j^{\pm}, j>1,\partial D^2=a_1^+$$
All boundaries are transversal to our system (see Fig 1). They
have numerical invariants $$\oint_{a_j}\omega^R=|a_j|>0, |a_1|\geq
|a_2|\geq \ldots \geq |a_g|>0$$
 The system does not have critical
points except $2g-2$ nondegenerate saddles. Every trajectory
starts at the {\bf in-boundary} $\bigcup a^+_j$ and ends at the
{\bf out-boundary} $\bigcup a_j^-$  (see Fig 4).

\begin{rem} Quite similar picture we obtain for the meromorphic 1-form
$\kappa$ on the Riemann 2-sphere $S^2=CP^1$ with 2g simple poles
(one of them at infinity),  and $2g-2$ simple zeroes. Such systems
$(\kappa)^R=0$ probably produce all topological types of
foliations in the domains like $\bar{V}$ obtained using
holomorphic forms on the Riemann surfaces with genus $g$.
\end{rem}

Consider a maximal abelian $Z^{2g}$-covering  $V'\rightarrow V$
 with basic shifts $a_j,b_j:V'\rightarrow V'$
for every class $T$ foliation. Its fundamental domain can be
obtained cutting $V$ along the Transversal Canonical Basis. The
connected pieces $A_j$ of the boundary exactly represent free
abelian groups $Z^2_j$ generated by the shifts $a_j,b_j$. At the
covering space a boundary component $A_j$ near this place looks
like standard square domain for 2-torus (see Fig 2). The foliation
near the boundary $A_j$ looks like a standard straight line flow
at the space $R^2_j$ with lattice $Z^2_j$. Our fundamental domain
$\bar{V}\subset V',\partial \bar{V}=\bigcup A_j$, is restricted to
the inner part of 2-parallelogram in every such plane $R^2_j$.
Topologically this domain  $\bar{V}$ is a 2-sphere with $g$ holes
(squares) , with
 boundaries $A_j$.

We can construct this covering analytically using  the Abel Map
$A=(A^1,\ldots,A^g)\in C^g$, with some initial point $P$:
$$A^j(Q)=\int_{P}^Q\omega_j,j=1,\ldots, g,P,Q\in V$$ Here
$\oint_{a_j}\omega_k=\delta_{kj}$ is a normalized basis of
holomorphic forms:
$$\oint_{a_j}\omega_k=\delta_{jk},\oint_{b_j}\omega_k=b_{jk}=b_{kj}$$.
The form $\omega=\sum u_k\omega_k$ is generic here. It defines a
one-valued function $F:V'\rightarrow
C$:$$F(Q')=A^{\omega}(Q')=\sum u_kA^k(Q')$$ where $Q'\rightarrow
Q$ under the projection $V'\rightarrow V$. For every component of
boundary of our fundamental domain we have for the basic shifts
$a_j,b_j:V'\rightarrow V'$:
$$F(a_j(Q'))=F(Q')+u_j,F(b_j(Q'))=F(Q')+\sum_k u_kb_{kj}$$ The
levels $$F^R=const$$ are exactly the leaves of our foliation
$\omega^R=0$ on the covering $V'$. The map $\pi_1(V)\rightarrow
H_1(V,Z)=Z^{2g}\rightarrow R$ is defined by the correspondence:
$$a_j\rightarrow u_k^R,b_j\rightarrow (\sum_ku_kb_{kj})^R$$ The
map $\psi: \pi^+_1(\Re)\rightarrow \pi_1(V)\rightarrow R^+$ of the
positive transversal semigroup (above) certainly belongs to the
semigroup $Z^{2g}_+$  where $(n,m)\in Z_+^{2g}$ if
$$\omega^R(\sum_{m_k,n_l}m_ka_k+n_lb_l)>0$$

There exists such choice of the phase vector $\eta_0$ that the of
$V'=A(V)\subset C^g$ satisfies to the equation identically:
$$\Theta (A(Q)-\eta_0|B)=0$$ for all points $Q$; the phase vector
$\eta_0$ depends on the initial point $P$ only. Here $$\Theta
(\eta^1,\ldots,\eta^g|B)=\sum_{n\in Z^g} \exp\{2\pi
i\sum_{k,j}b_{jk}n_kn_j+\sum_jn_j\eta^j\}$$ For $g=2$ it is a
complete equation defining this submanifold as a {\bf
$\Theta$-divisor} in the Jacobian variety.

\vspace{1cm}

\centerline{\Large Section 2. Some General Statements.}

\vspace{0.5cm}

Consider now any compact nonsingular algebraic curve $V$ with
holomorphic generic 1-form $\omega=\sum u_k\omega_k$ and foliation
$\omega^R=0$. Our foliation is defined through the complex
analytic function on the abelian covering $F:V'\rightarrow C$,
where this function is defined by the integral along the path
joining the initial point $P$ with a variable-point $Q'$ which is
a pair $(Q\in V,[\gamma])$ where $\gamma$  is a homology class of
paths joining $P,Q$: $$F(Q')=\int_{\gamma}\omega,Q'\rightarrow Q$$
This is a restriction of the linear function $\sum u_kA^k$ in the
space $C^g$ generated by the normalized basis of holomorphic forms
$\oint_{a_q}\omega_j=\delta_{qj}$, to the complex curve
$V'=A(V)\subset C^g$. Here $A$ is a multivalued Abel Map. The
levels $F^R=const\in R $ are exactly our covering trajectories in
$V'$. The critical values are $F(S)\subset C$ where $dF|_{Q'}=0$
for $Q'\in S$.

{\bf Definitions}.

1.We call foliation {\bf Generic}  if it satisfies to the
following requirements: It has only nondegenerate critical points
(i.e. saddles);  There exist no saddle connections  (no
separatrices  joining two saddles); Even more, no one line of our
selected parallel family  $F=const$ in $C$ crosses the critical
value set $F(S)$ twice: no one line of this family crosses twice
also the ''quasilattice''  $Z^{2g}\subset C$ generated by $2g$
complex numbers in $C=R^2$ where
$u_k,\sum_ku_kb_{kl},k,l=1,\ldots,g$; Every periodic trajectory is
homologous to zero (i.e. it divides Riemann Surface into 2
pieces). Without any further quotations we are going to consider
only {\bf The Irreducible Generic Foliations} which do not have
periodic orbits at all.

2.By the {\bf Almost Transversal Curve} we call every parametrized
piecewise smooth curve consisting of the two type smooth pieces:

First Type: Moving transversally to foliation in the same
direction.

Second Type: Moving along the trajectories of foliation in any
direction.

A simple lemma known many years claims that {\bf every almost
transversal curve can be approximated by the smooth transversal
curve} with the same endpoints (if there are any). In many cases
below we construct closed almost transversal curves and say
without further comments that we constructed smooth transversal
curve.

3. We call by the {\bf Plane Diagram} of foliation of the type
$T^k$ with transversal canonical basis a {\bf Topological Type} of
foliation on the Riemann Surface $\bar{V}$ obtained from $V$ by
cuts along this basis (see Fig 10 for $g=2$, and the descriptions
below).

\newtheorem{thm}{Theorem}

\begin{thm} Every generic foliation given by the holomorphic 1-form $\omega^R=0$
on the algebraic curve of genus 2 belongs to the class $T$, i.e.
admits a full Transversal Canonical Basis
\end{thm}

Proof of this theorem follows from the following two lemmas:

\begin{lem}Every generic foliation $\omega^R=0$ belongs to the
 class $T^0$ for genus equal to  2
\end{lem}

Proof. Take any trajectory such that its limiting set in both
directions  contains at least one nonsingular point. In fact,
every nonseparatrix trajectory has this property: if its limiting
set contains critical point, it also contains a pair of
separatrices entering and leaving it. Nearby of the limiting
nonsingular point our trajectory appears infinite number of times.
Take two such nearest returns and join them by the small
transversal segment. Obviously, this closed curve consisting of
the piece of trajectory  and small transversal segment, is an
almost transversal curve: it can be approximated by the closed
non-selfintersecting smooth transversal curve. We take this curve
as a cycle $a_1$. Now we cut $V$ along this curve and get the
surface $\tilde{V}$ with 2 boundaries
$\partial\tilde{V}=a_1^+\bigcup a_1^-$. We take any trajectory
started at the cycle $a_1^+$ and ended at $a^-_1$. Such trajectory
certainly exists. Join the ends  of this trajectory on the cycle
$a^-_1$ by the positive
 transversal segment  along the cycle $a_1$ in $V$. We
get a transversal non-selfintersecting cycle $b_1$, crossing $a_1$
transversally in one point. Cut now $V$ along the pair $a_1,b_1$.
We get a square $\partial D^2=a_1b^{-1}_1a_1^{-1}b_1\subset R^2$
with  a 1-handle attached to the disk $D^2$ inside. We choose
notations for cycles in such a way that our foliation enters this
square along the piece $A_1^+=a_1b^{-1}_1\subset
\partial D^2$, and leaves it along the piece
$A_1^-=a_1^{-1}b_1$ (see Fig 3). Nearby of the angles where these
pieces $R^{\pm}$ are attached to each other, our trajectories
spend small time inside of the square entering and leaving it (see
Fig 3). So moving inside from the both ends of the segment
$A_1^+$,, we find 2 points $x^+_{1,1},x^+_{2,}\in A_1^+$ where
this picture ends (because the genus is more than 1): These points
are the ends of the separatrices of the saddles. Very simple
qualitative arguments show that $x^+_{1,1}\neq x^+_{2,1}$.  Take
any point between $x^+_{1,1}$ and $x^+_{2,}$ (nearby of
$x^+_{1,1}$). The trajectory started at this point crosses $A_1^-$
somewhere (see Fig 3). Join the end-point $y'$ of this piece of
trajectory by the transversal segment $\sigma=y'y$ along the curve
$A_1^-$ in positive direction with the point equivalent to the
initial one. This is a closed curve $c_1$ transversal to our
foliation. If initial point is located on the cycle $a_1$, the
curve $c_1$ does not cross the cycle $b_1$. We take cycles
$b_1,c_1$ as a basis of the transversal $a$-cycles. If initial
point is located on the cycle $b_1$, the curve $c_1$ does not
cross $a_1$. In this case we take $a_1,c_1$ as a basis of the
transversal $a$-cycles. It is easy to see that $c_1$ cannot be
homologous to $b_1$ in the first case. Therefore it is a right
basis of the $a$-cycles which is transversal. The second case is
completely analogous. Our lemma is proved.

\begin{lem} For every generic foliation $\omega^R=0$ of the class $T^0$
on the algebraic curve of genus 2,
 the transversal basis of $a$-cycles can be extended to the full
Transversal Canonical Basis $a,b$, so every $T^0$-class foliation
belongs to the class $T$.
\end{lem}

Proof. Cut the Riemann Surface $V$ along the transversal cycles
$a_1,a_2$. Assuming that $\oint_{a_1}\omega^R=|a_1|\geq
|a_2|=\oint_{a_2}\omega^R$, we realize this domain as a plane
domain $D^2_*$ as above (see Fig 4 and the previous section). Here
an external boundary $\partial_{ext} D^2_*=a_1^+$ is taken as our
maximal cycle. The elementary qualitative intuition shows that
there are only two different topological types of the plane
diagrams (see Fig 4, a and 4,b): The {\bf first case} is
characterized by the property that for each saddle all its
separatrices end up in the four different components of boundary
$a_1^{\pm},a_2^{\pm}$. We have following matrix of the trajectory
connections of the type $(k,l):a^+_k\rightarrow a^-_l$  for the
in- and out-cycles and their transversal measures:

 $a^+_1\rightarrow a_1^-$ (with measure $a$), $a^+_k\rightarrow a^-_l, k\neq l$
(with measure $b$ for $(k,l)=(1,2),(l,k)=(2,1)$), and
$a^+_2\rightarrow a^-_2$ with measure $c$. All measures here are
positive. We have for the measures of cycles:
$|a_1|=a+b,|a_2|=c+b$. This topological type does not have any
degeneracy for $a=c$.

The diagonal trajectory connections of the type $(l,l)$ generate
the transversal $b$-cycles closing them by the transversal pieces
along the end-cycles in the positive direction.

In the {\bf second case} we have following matrix of trajectory
connections:

$a_1^+\rightarrow a^-_1$ with measure $a>0$, $a^+_1\rightarrow
a^+_2$ with measure $b>0$, $a_2^+\rightarrow a_1^-$ with the same
measure $b>0$. We have $b=|a_2|,a+b=|a_1|$. So we do have the
trajectory-connection $a_1^+\rightarrow a_1^-$, but we do not have
the second one, of the type $(2,2)$. However, we may connect
$a_2^+$ with $a_2^-$ by the almost transversal curve as it is
shown in the Fig 4,b), black line $\gamma$.  So we construct the
cycles $b_1,b_2$ as the transversal curves crossing the cycles
$a_1,a_2$ only.

Our lemma is proved.

 Therefore the theorem is also proved.

{\bf The case $g=3$}. A lot of concrete foliations of the class
$T^2$ were demonstrated above for the real nonsingular algebraic
curves $$w^2=\prod (z-z_1)\ldots (z-z_8)=R(z), z_j\in R$$, with
the cycles $$a_1=[z_2z_3], a_2=[z_7z_8], a_3=[z_4z_5],
b_1=[z_1z_2], b_2=[z_6z_7]$$ and $\omega=P_2(z)dz/R(z)^{1/2}$ The
polynomial $P_2(z)=u+iv$ should be chosen such that its real part
does not have zeroes in the segments $[z_2z_3], [z_4z_5],
[z_6z_7]$, and its imaginary part does not have zeroes in the
segments $[z_1z_2],[z_7z_8]$ (the proof is identical to one in the
Section 1 for the constant $u,v$).

{\bf Question:} Is it possible to extend this basis to the
Transversal Canonical Basis? As we shall see below, the answer is
negative in some cases: we need to reconstruct our incomplete
basis in order to extend it to the full transversal canonical
basis.

 Let a nonsingular algebraic curve $V$ of the genus $g=3$ is
given with the generic foliation $\omega^R=0$ belonging to the
class $T^2$ with the incomplete  basis $a_1,a_2,a_3,b_1,b_2$
transversal to the foliation. We construct its {\bf Plane
Diagram}. After cutting the Riemann Surface along these cycles we
realize it as a plane domain with following components of the
boundary:

 The external boundary $$A_1= (a_1b_1^{-1})\bigcup
(a_1^{-1}b_1)=A_1^+\bigcup A_1^-$$

 The internal boundary
$$A_2=(a_2b^{-1}_2)\bigcup (a_2^{-1}b_2)=A_2^+\bigcup A^-_2$$

  The
interim boundary $a_3^+\bigcup a_3^-$ inside.

 Our notations are chosen in such a
way that trajectories  enter our domain through the piece with the
sign $+$ and leave it through the pieces with the sign $-$. Our
Hamiltonian provides a transversal measure. We make a numeration
such that: $$2|A_1^{\pm}|=A_1> A_2=2|A_2^{\pm}|,
|a_3|=|a_3^{\pm}|=a$$ Here $A_k$ means also the measure of this
boundary component.  Nearby of the ends of the segments $A_l^+$
the trajectories enter our domain and almost immediately leave it
through the piece $A_l^-$. Therefore, there exist the first points
in $A^{\pm}_l$ where this picture ends. These are the endpoints
$x_{1,j}^{\pm},x_{2,j}^{\pm}\in A_j^{\pm}$ of the pair of
separatrices of  saddles (see Fig 5). We don't see here the other
pair of separatrices for these saddles.

\begin{lem}
A complete list of  topologically  different types of the Plane
Diagrams in the class $T^2$ for the genus $g=3$ can be presented.
It shows that there is only one type such that  we cannot extend
the incomplete transversal basis to the complete transversal basis
(see Fig 6,a)): no closed transversal curve exists in this case
crossing the cycle $a_3$ and not crossing other curves
$a_1,b_1,a_2,b_2$ (i.e. joining $a^+$ and $a^-$ on the plane
diagram). In all other cases such transversal curve $b_3$ can be
constructed.
\end{lem}

 Consider now a special case of the class $T^2$ where the
transversal incomplete basis cannot be extended (Fig 6,a).
\begin{lem}
There exists a  reconstruction of this basis such that the new
basis can be extended to the complete transversal canonical basis

\end{lem}

For the proof of the second lemma, we construct
 a closed transversal curve $\gamma$  such that
it crosses the cycle $a=a_3$ in one point and  crosses  the
segment $A^+$  leaving our plane diagram (see Fig 7). It enters
$A^-$ in the equivalent point, say, through the cycle $b_1$. We
take a new incomplete transversal basis $a_2,b_2, a_3, \gamma,
a_1$. If $\gamma$ crosses $A^{\pm}_1$ through the cycle $a_1$, we
replace $a_1$ by the cycle $b_1$ as a last cycle in the new
incomplete basis of the type $T^2$. The proof follows from the
plane diagram of the new incomplete basis (we drop here these
technical details, especially the list of the plane diagrams
implying lemma 5 ).

 Comparing these lemmas with the construction
of special foliations in the previous section on the real Riemann
Surface, we are coming to the following

{\bf Conclusion.} For every real hyperelliptic Riemann Surface of
the form $w^2=\prod (z-z_1)\ldots (z-z_8)=R(z), z_k\neq z_l\in R$,
every generic form $\omega= P_2(z)dz/R(z)$ defines foliation
$\omega^R=0$ of the class $T$ if real and imaginary parts of the
polynomial $P=u+iv$ do not have zeroes on the cycles indicated
above in the lemma 2 for   genus $g=3$.

 We prove in the Appendix that every generic foliation admit a
 Transversal Canonical Basis for $g=3$. G.Levit communicated to
 the author another proof of this theorem for all $g\geq 2$ based on
 the
 construction of  pant decomposition of  Riemann surface with
 $3g-3$
 boundary curves transversal to dynamical system
 (his proof is also included  in the Appendix).
 \vspace{0.3cm}

\centerline{\Large Section 3. Topological Study of the Class $T$.}

\vspace{0.5cm}

Let us describe here some simple general topological properties of
the class $T$ foliations $\Re$ and corresponding Hamiltonian
Systems. After cutting the surface $V$ along the transversal
canonical basis $a_j,b_j$, we are coming to the fundamental domain
$\bar{V}$ of the group $Z^{2g}$ acting on the maximal abelian
covering $V'\subset C^g$ imbedded by the Abel Map. It leads to the
Plane Diagram $D^2_*$ with $g$ boundary ''squares'' $\partial
D^2_*=\bigcup_{j} A_j$ where $A_j= a_jb^{-1}_ja_j^{-1}b_j
=A^+_j\bigcup A^-_j$. Every piece $A_j^{\pm}$ consists of  exactly
two basic cycles $a_j,b_j$ attached to each other. These pieces
are chosen such that trajectories enter $A^+_j$ from outside
through $a_jb_j^{-1}$ and leave it into the fundamental domain
$\bar{V}$  except the areas nearby of the ends. The trajectories
enter $A^-_j$ from inside and leave fundamental domain. There is
exactly $2g-2$ saddle points inside of $\bar{V}$. They are not
located on the selected transversal cycles $a_j,b_j$. Our
foliation nearby of each boundary square $A^+$ looks exactly as a
straight line flow. It means in particular that there exist two
pairs of points $x_{1,j}^{\pm},x^{\pm}_{2,j}\subset A^{\pm}_j$
which are the endpoints of separatrices in $\bar{V}$, nearest to
the ends at the each side $A^{\pm}_j$ (see Fig 3 and 5). We call
them {\bf The Boundary Separatrices} belonging to {\bf The
Boundary Saddles } $S_{j,1},S_{j,2}$ for the cycle $A_j$. We see
$2g$  of such saddles $S_{j,1},S_{j,2}$ looking from the
boundaries of $A_j$,  but some of them are in fact the same. At
least two of them should coincide leading to the saddles of the
types $<jjkk>$ where all incoming and outcoming separatrices are
of the boundary type (The ''Double-boundary'' saddles).

{\bf Definitions} 1. We say that the saddle point $S \in \bar{V}$
has a type $<jklm>$ if it has two incoming separatrices starting
in   $A^+_j,A^+_l$ and two outcoming separatrices ending in
$A^-_k,A^-_m$. The indices are written here in the cyclic order
corresponding to the orientation of $\bar{V}$. Any cyclic
permutation of indices defines an equivalent type. we normally
write  indices of the type  starting from the incoming separatrix,
as $<jklm>$ or $<lmjk>$.

 2. We call
foliation {\bf Minimal} if all saddle points in $\bar{V}$ are
contained in the set  of boundary saddles. In particular, their
types are $<jjkl>$. We call foliation {\bf  Simple} if there are
exactly two saddle points of the double-boundary  types $jjkk$ and
$jjll$ correspondingly. The index $j$ we call {\bf Selected}. We
say that foliation has a rank equal to $r$, if there exists
exactly $r$ saddles of the types $<jklm>$ where all four indices
are nonboundary. In particular, we have $0\leq r\leq g-2$. There
is exactly $ t$ saddles of the double-boundary types like $jjkk>$
where $t-r=2$. There is also $2g-2t$ other saddles of the types
like $<jjkl>$ where only index $j$ corresponds to the boundary
separatrices. The extreme cases are $r=0,t=2$ which we call {\bf
special} (above), and $r=g-2, t=g$ which we call {\bf maximal}. In
the maximal case there exists a maximal number of saddles whose
separatrices arrived from the nonboundary parts of $A^+_k$, and
all boundary type saddles are organized in the pairs. One might
say that for the maximal type every index is selected. {\bf For
maximal type the genus should be an even number} because the
boundaries $A^j$ are organized in the cycles now, and every cycle
should contain even number of them, by the elementary orientation
argument. The relation $2g-t+r=2g-2$ for the total number of
saddles gives $t-r=2$. For the case $g=2$ we obviously have $r=0$.
For the case $g=3$ the only possible case is $t=2,r=0$ (the
simple foliations); the case $t=3,r=1$  cannot be realized for
$g=3$ because it is maximal in this case: however, the  maximal
case corresponds to the even genus only. Therefore it is available
only for the genus not less that 4.

\vspace{0.2cm}

{\bf How to build these   systems topologically?}

In order to answer this question, let us introduce  following {\bf
Building Data} (see Fig 9):

I.{\bf The Plane Diagram} consisting of  the generic Hamiltonian
System on the 2-sphere $S^2$ generated by the hamiltonian $H$ with
nondegenerate critical points only (centers and saddles) sitting
on different levels. Let one center is located in the point $0$,
and another one in $\infty$. It has $t$ centers and $r$ saddles.
Let exactly $g$ transversal oriented segments are given
$t_1,...,t_g \subset S^2$ with transversal measures
$m_1,m_2,...,m_g$ provided by hamiltonian, such that:

a.They do not cross each other;  the values of Hamiltonian in
their ends, centers and saddles  are distinct except that exactly
 two of them
meet each other in every center; They do not touch any saddle
point on the two-sphere.

b.Every cyclic and separatrix trajectory of the hamiltonian system
on $S^2$
 meets at least one of these segments.

 We make cuts along these segments and define the sides $t_j^{\pm}$
 where trajectories  leave and enter it correspondingly.

II.{\bf The Torical  Data} consisting of the $g$ tori
$T^2_j,j=1,...,g$ with distinct hamiltonian irrational straight
line flows
 and selected oriented transversal
segments $s_1,...,s_g$ (one for each torus). Their  transversal
measures are equal to the same numbers $m_1,...,m_g$. The
Transversal Canonical Basis  $a_k,b_k$ in every torus is selected
where $a_k$ are positive, $b_k$ are negative, and
$|a_k|+|b_k|>m_k$.

We make similar cuts along these segments in the tori, and define
their sides $s^{\pm}$ in the same way.

Identifying the segments $s^+_j$ on the tori with $t^-_j$ on the
plane and vice versa, we obtain a Riemann surface $M^2_g$ with
foliation which has a transversal measure. The explanation should
be given concerning the centers and the ends of the segments:

We construct our gluing in such a way that every end of the
segment $t_j\subset S^2$ defines exactly one saddle of the
boundary type $<jjkl>$. Here $t_k,t_l$ are  the  segments  joined
by
 the pieces of the same  trajectory with  the end of the segment $t_j$ on the
2-sphere $S^2$. These pieces of  trajectory provide  a pair of
nonboundary separatrices for the saddle on the Riemann Surface. We
assume that they meet these segments in the inner points because
the  foliation is generic.

By definition, every center  generates
 a double-boundary
saddle of the type $<jjkk>$. So we have $t-r=2$.

In order to obtain this set of data  from the generic foliation
given by the real part of holomorphic one-form with transversal
canonical basis, we perform following operations:

1.Cut our surface along the TCB. The boundary of this domain
$\bar{V}$  is equal to the union $\bigcup_j (A^+_j\bigcup
A^-_j)=\partial \bar{V}$. Every component is presented as lying in
the Fundamental Parallelogram $P_j$ of the 2-torus $T^2_j$. Our
flow covers the boundary  of $P_j$ as a straight line flow: the
trajectories enter through the path realized by the pair of cuts
$a_jb_j^{-1}$ and leave through the path $a_j^{-1}b^j$ (see Fig
2).

2.Find for every 2-torus $T^2_j$ the pair of  boundary saddles in
$P_j$ and join them by the pair of transversal segments $s_j^{\pm}
$. They should meet each other  in the boundary saddles only (see
Fig 5). We perform this operation in the fundamental parallelogram
$P_j$ representing our torus. These segments should be chosen in
such a way that outside of them in $P_j$ (or in the plane
$C=\bigcup_{g\in Z^2}g(P_j)$ near the one-skeleton), we have a
straight line flows.

3.Cut our surface $M^2_g$ along all these segments $s^+_j\bigcup
s^-_j$.  It is divided now into the torical pieces $T^2_j$ and one
plane piece $S^2$ whose boundary consists of the curves
  $\bigcup_j(t^+_j\bigcup t^-_j)$  for $S^2$, and $s^+_j\bigcup s^-_j$
 for the tori $T^2_j$. After cutting the surface $M^2_g$ along the
 pieces $s_j^{\pm}$ we keep the notation $s^{\pm}_j$ for the curves
 in the tori, but for the plane part $S^2$ we change notations for
  these curves, and
 denote them by $t^{\pm}_j$.

4 Now we glue $t^+_j$ with $t^-_j$ for the sphere $S^2_j$,  and
$s^+_j$ with $s^-_j$ for the tori $T^2_j$, preserving the
transversal measure.  The system on the
 2-sphere appears with $g$ selected transversal segments $t_j$.
We have also $g$ 2-tori $T^2_j$ with the straight line flows and
transversal pieces $s_j$ whose measures are equal  $m_j$.

5.Near the double-boundary saddles we are coming to the picture
topologically equivalent to the center, but this equivalence in
non-smooth.

 We can see  that our construction allows to imitate all
topology of foliation. It  preserves also the measure-type
invariants.

Therefore we are coming to the following

\begin{thm} Every generic foliation given by the real part of
 holomorphic one-form,
can be obtained by the measure-preserving gluing of the
 pieces $(S^2,H,t_1,...t_g)$
and $(T^2_1,s_1),...,(T^2_g,s_g)$ along the transversal segments
$t_j$ and $s_j$, as it was described above. For the genus $g=2$ we
cam remove  sphere $S^2$ from the description: Every generic
foliation given by the real part of holomorphic one-form, can be
obtained from the pair of tori $(T^2_1,s_1)$ and $(T^2_2,s_2)$
with different irrational straight line flows and  transversal
 segments $s_1,s_2$ with transversal measure $m_1=m_2=m$
\end{thm}

{\bf Example. The Topological Types of the Minimal Foliations}.

 For
the Minimal Foliations above we have $t=2,r=0$. The hamiltonian
system on the 2-sphere is trivial (see Fig 9,a): It can be
realized by the rotations around the point $0$. There are no
saddles on the sphere here, and the second center we take the
point $\infty$. For the simplest case $g=2$ there are two segments
$t_1,t_2$. Both of them join $0$ and $\infty$. So they form a
cycle of the length 2. The difference between the values of
Hamiltonian $H(\infty)-H(0)$ is equal to $m_1=m_2$. All possible
pictures of the transversal segments can be easily classified here
for every genus $g$ (see Fig 9,a and b).

There exist following types of topological configurations only:

a.The Plane Diagram has exactly one ''cycle'' $t_1,t_2$ of the
length two (like for $g=2$) and $g-2$ disjoint segments $t_j,j\geq
3$; This type is available for all $g\geq 2$.

b.The Plane Diagram has two pairs $t_1,t_2$ and $t_3,t_4$ where
the members of each pair meet each other either  in the center $0$
or in $\infty$, and $g-4$ disjoint segments $t_j,j\geq 5$. In the
second case we have $g\geq 4$.

c.The Plane Diagram has exactly one connected set consisting of 3
segments passing through both centers $t_1t_2t_3$ and $g-3$
disjoint segment $t_j,j\geq 4$. For this type we have $g\geq 3$.

\begin{thm} For $g=2,3$ every class $T$ foliation is simple. A maximal
type exists only for even genus $g\geq 4$
\end{thm}

Proof. For $g=2$ this is obvious and was already established
above: all generic (irreducible) hamiltonian foliations are
simple. Both of indices $k=1,2$ are selected; at the same time
they are maximal. Consider now the case $g=3$. Our foliation can
be either simple ($t=2,r=0$) or maximal ($t=3,r=1$) in this case.
We have $t=3,2g-2t=0$ for the  maximal case. If it is so, every
boundary saddle should be paired with some other. So there is a
cyclic sequence of boundary saddles containing all three boundary
components. However, every cyclic sequence
 should contain even number of boundary components ( and the same
number of boundary saddles), otherwise the orientation of
foliation is destroyed. This is possible only for even number of
indices which is equal to genus. Our conclusion is that $g\geq 4$.
This theorem is proved.

The 2-sphere is covered by the nonextendable {\bf '' corridors''}
between two transversal segments $t_j,t_k\subset S^2$. They are
the strips of nonseparatrix trajectories  moving from the inner
points of $t_j$ to the inner points of $t_k$ not touching  any
points of $t_l$  and the saddles. The right and left sides of
these corridors are either separatrices of the saddles in $S_2$ or
the trajectories passing on $S^2$ through the ends of some
segments $t_l$.

Classification of the generic Morse functions $H$ on the sphere
$S^2$ can be given easily (see Fig 9, c): Take any connected
trivalent finite tree $R$. It has vertices divided into the  $r$
{\bf Inner Vertices} and $t$ {\bf Ends}. Assign to each vertex
$Q\in R$ a value $H(Q)$; these values   are not equal to each
other $H(Q)\neq H(P),P\neq Q$; therefore the edges become
oriented,  looking ''up'', to the direction of increasing of $H$.
Every inner vertex $Q$ is a ''saddle'', i.e.
 $$\min_iH(Q_i)<H(Q)<
\max_iH(Q_i), i=1,2,3$$ where $Q_i$ are the neighbors. The
function $H$ on the graph should be such that for every edge
$[Q_1Q_2]\subset R$  we have $$H(Q_1)\leq H(P)\leq H(Q_2)$$ where
$H(P)$ is monotonic in this edge. We may take it linear.

For the description of the set of transversal segments $t_j\subset
S^2$ we introduce a {\bf locally constant set-valued function}
$\Psi (P),P\in R $ on the graph $R$: the values $\Psi(P)$ are the
cyclically ordered nonempty finite sets
$$\Psi(P)=\{t_1(P)<t_2(P)<...<t_q(P)<t_1(P)\}$$ such that: The
number $q=|\Psi(P)|$ of the points $t_j(P)$ in this set, is equal
to 2 nearby of the ends $Q_j$ collapsing to one in the endpoint
$t_1(P)=t_2(P)$ if $P=Q_j$. The number $q$ may change
$q\rightarrow q\pm 1$ in the isolated points $P_l\in R$ such that
$H(P_l)\neq H(Q)$ for all inner vertices $Q$. Passing ''up''
through the inner vertex $Q$, this set either splits into the pair
of cyclically ordered sets $\Psi(P)\rightarrow \Psi_1\bigcup
\Psi_2$ inheriting orders  where $|\Psi_1|+|\Psi_2|=|\Psi(P)|$, or
some pair of cyclically ordered sets is unified into the one
ordered set choosing some initial points in each of them (the
inverse process). Every continuous ''one-point branch'' $t_j(P)$
living in the  vertical path between the points $P_1,P_2\in R$,
defines the transversal segment $t_j$. Its transversal measure is
equal to $m_j=H(P_2)-H(P_1)>0$

One may imagine that the graph $R$ is imbedded into the space
$R^3$. The sphere $S^2$ appears as a boundary of the small
$\epsilon$-neighborhood of $R$ in $R^3$. The function $H$ should
be realized as a ''height function'' $S^2\rightarrow R$. The
points $t_j(P)\in\Psi(P)$ are marked on the boundary of this small
neighborhood. The nearest component of the level $H=c$ is exactly
a small circle near the noncritical point $P\in R$. Changing $c$,
the points $t_j(P)$ are varying according to the rules above.

{\bf Every trajectory $\gamma\in V$ defines a sequence of elements
$$...W_{-M}W_{-M+1}....W_NW_{N+1}...=W(\gamma)$$ where $W_N\in
H_1(V,S^2)=H_1(V)$} We compute these elements $W_N$ through the
3-street model and $m_j$-dependent new TCB in the next chapters.

\vspace{0.3cm}

\centerline{\Large Section 4. The Three-Street Picture on the
Torus } \centerline{\Large The Case $g=2$. The Maximal Case for
$g=4$ }

\vspace{0.5cm}

Let us describe the  case $g=2$  more carefully. We have here two
selected levels $1,2$. The reduced or working measure of the
transition from $C_1$ to $C_2$ and back is equal to $m$. The
transition map is an orientation preserving isometry of segments
$$\Phi:s=s_2\rightarrow s_1=s$$ We have two planes $C_1,C_2$ with
two lattices $Z^2_1,Z^2_2$. A family of parallelograms starts in
the selected points. The transversal segments $s=s_2,s=s_1$ are
located in each of them. Everything is repeated periodically in
each space $C_k$ with its own lattice. Our data include
$|a_j|,|b_j|$ for $j=1,2$ and the transition measure $m$.

Consider the vertical flow in $C_k$.

{\bf Question}. How long the trajectory can move in $C_k$ (i.e. in
the torus $T^2_k$) until it hits some periodically repeated copy
of the  segment $s$?

It starts and ends in some segments of the selected periodic
family generated by the segment $s$  in the corresponding
parallelograms not crossing any segments in between. Such paths
with fixed ends  form the connected strip. We require that these
strips cannot be extended to the left and right: every trajectory
in the strip ends in the same  segments. The extension of the
strip to the right or to the left   meets some saddles. They are
presented by the ends of the segments of our family. Every such
strip has a {\bf Height $h$} and {\bf Width $w$}. The width is
equal to the transversal measure. The height depends on the
lattice periods. It  has a meaning only as a topological quantity
$h\in H_1(V,Z)$.

 {\bf Definition}.We call the unextendable strips by {\bf The Streets}
 and denote them $p^{\tau}_k, k=1,2$. We  denote a longest unextendable strip
by  $p^0_k$. It   meets the ends of some segments of our family
strictly inside along the segment number zero (see Fig 11). The
upper and lower segments of this strip should be located also
strictly inside of the corresponding segments $s, s'''$. Their
Heights and Widths we denote by $h^{\tau}_k, |p^{\tau}_k|$
correspondingly. The street number 2 is located from the right
side from the longest one, the street number 1--from the left
side.

\begin{lem}
For every foliation with transversal canonical basis and for both
planes  $k=1,2$ there exists exactly three streets
$p_k^{\tau},\tau=0,1,2$ such that $\sum_{\tau } |p^{\tau}_k|=m$
and $h^1_k+h^2_k=h^0_k$. Two smaller streets are attached to the
longest one from the right and left sides. This picture is
invariant under the involution changing time and orientation of
the transversal segments. The union of the three streets started
in the segment $s$, is a fundamental domain of the
 group $Z^2_k$ in the plane $C_k$ (see Fig 11).

\end{lem}

\begin{rem} Another fundamental domain associated with segment $s$
can be constructed in the  form of the ''most thin
parallelogram''. It is generated by two vectors depending on  the
measure of this segment. The first vector corresponds to the shift
$s\rightarrow s'$ where $s'$ is the second end of the street
number 1 started in $s$ (from the left sight of the longest street
number 0--see Fig 11). The second vector corresponds to the shift
$s\rightarrow s''$ where $s''$ is the second end of the street
number 2
 (from the right side of the street number  0--see Fig 11).
We shall discuss this parallelogram later, in the last .section.
This $m$-dependent ''thin'' basis of the lattice $Z^2$ can be
canonically lifted to the free group with 2 generators $a,b$. We
denote them $a^m(+)=a^*,b^m(+)=b^*$ such that the product path
$a^*b^*(a^*)^{-1}(b^*)^{-1}$
 contains exactly one segment $s$ inside. This new generators are
 also transversal to our foliation, so we can construct a new TCB out of them.
 We shall use them instead of the original TCB
 because they are adjusted to the idea of our description of
foliations with TCB on the surfaces of higher genus $g>1$. Their
transversal measures are following: $$|a^m_k|=|p^0_k|+|p^2_k|,
|b^m_k|=|p^0_k|+|p^1_k|$$ where
$\sum_{\tau}|p^{\tau}_k|=m,\tau=0,1,2,k=1,2 $

\end{rem}

Proof. Construct first the longest street $p^0_k$. We start from
any nonseparatrix trajectory ending in some segments inside both
of them. Extending this strip in both (left and right) directions,
we either meet the ends of upper or lower segments or meet some
segment whose height is strictly between. In the first case we see
that after passing the most left end we can construct longer
trajectories. Do it and start the same process with longer strip.
Finally, we reach the locally maximal vertical length of the
height. After that we extend it to the right and left. We
necessarily meet from both sides some segments with heights in
between, otherwise its height cannot be locally maximal. So the
maximal provider is constructed. Consider the neighboring streets
from the right and left sides. This is exactly two other streets.
Their heights are smaller. Denote the left one by $p^1_k$ and the
right one by $p^2_k$. What is important, is that the neighboring
street can be extended till the left end of the upper segment of
the locally longest one. This statement follows from the
periodicity of the system of segments: we cannot meet any segment
from the left until we reached the end of the upper one; the newly
met segment should have the locally longest one from its right
side. The same argument can be applied to the right extension and
to the lower parts as well. We can see that our locally longest
street is surrounded by the exactly  four unextendable domains
(two domains from each side). They are restricted from the right
and left sides by the ends of the lower or upper segments (see Fig
11 ). The pairs of domains located across the diagonal of each
other are equal. All our relations immediately follow from that.
There are no other unextendable streets except these three.

Lemma is proved.

{\bf The Model of the Motion for $g=2$} is following: The motion
is vertical (up). Take two copies of the horizontal segment
$s=s_k$ of the length $m$. Put each of them  into the plane $C_k$
for $k=1,2$ inside of the parallelogram $P_k$. Construct three
vertical streets $p^{\alpha}_k,k=1,2$ over $s=s_k$ from the upper
side (the longest street lies in between). Assume that the segment
$s=s_k$ belongs to the parallelogram $P_k$ with $Z^2$-index
$(0,0)$ in both cases. Assign to the upper end of each street
$p^{\alpha}_k$ an integer 2-vector $h^{\alpha}_k\in Z^2$. It is
exactly a $Z^2$- index of the lattice parallelogram where it is
located. We have $h^1_k+h^2_k=h^0_k$. The trajectory starts in
$C_k$ at the segment $s_k$; it moves along the street vertically.
After reaching the upper end, it jumps to $C_l,l\neq k$, exactly
to the same point of the lower segment (modulo periods) as the
endpoint in the street (but on the different plane). After that it
moves along the corresponding street $p^{\beta}_l$ in the new
plane $C_l$, and so on. Remember that the widths (i.e. the
transversal measures) of the
 $l$-streets are different satisfying only to the conservation law
that their sum is also equal to $m$.  After each period of the
straight-line  motion along the street $p^{\alpha}_k$, the
$Z^2_k$-index of fundamental domain $\bar{V}$ changes: we add a
vector $h^{\alpha}_k$ to the $Z^{2n}$-number of our domain
$\bar{V}$ changing only  the pair of components corresponding to
the plane $C_k$. Other components remain unchanged. During the
jump $\Phi$ from $C_1$ to $C_2$ and back all $Z^{2n}$-numbers of
fundamental domain $\bar{V}$ remain unchanged.

This leads to the full topological description of foliation in the
combinatorial form for $g=2$ if we can calculate all incoming
homotopy and homology classes of streets effectively.

 {\bf How to describe time
dynamics? How much time is needed to pass every street?} For the
precise definition of time intervals we mark a pair of segments
$s_2 ,s_1$ transversal to foliation (except the ends) leading from
one saddle to another in the fundamental domain $\bar{V}$ . By
definition, what is presented by the ''streets'' combinatorially
in $C_k$, is presented in the actual dynamical system  by the
strips of trajectories which start in $s_k$  and end up in
$s_l,l\neq k$. The boundary separatrices of shorter streets
$p^{\alpha}_k,\alpha=1,2,$ have exactly one saddle points at each
of their boundaries in the lower or upper ends of the street. The
longest streets $p^0_k$ have also one saddle at each of their
boundaries  located somewhere inside of the ends (see Fig 11).

We define {\bf The Time Characteristic Functions} $t(x)>0$ for
$x\in (0,|p^{\alpha}_k|)$ for passing every street where
$t\rightarrow +\infty$ if $x\rightarrow 0$ or $x\rightarrow
|p^{\alpha}_k|$. In all cases the asymptotics is like $t(x)\sim
-c\ln x$ or $t(p^{\alpha}_k-x)\sim -c\ln x$, but the constants are
different. Every saddle provides two positive constants
$c_m>0,m=1,2$, where $m=1$ corresponds to the left side, and $m=2$
to the right side of the street. We have following asymptotics,
$x\rightarrow +0$

1.For $x\in p^0_k$ $$t(x)\sim (-c_2) \ln x; t(p^{0}_k-x)\sim
(-c_1) \ln x$$

 2.For $x\in p^{1}_k$: $$t(x)\sim (-c_2/2)\ln x;
t(p^{1}_k-x)\sim (-c_{1}/2)\ln x$$

3.For $x\in p^2_k$ $$t(x)\sim (-c_1/2)\ln x; t(p^2_k-x)\sim
(-c_2/2)\ln x$$

The exact value of the time characteristic functions $t(x)$ should
be calculated numerically. Their singularities are important, for
example,  for the ergodic properties of hamiltonian systems.  How
much time the trajectory spends in this or that area?   What size
fluctuations might have?  For $g=1$ in the presence of saddles
this problem was studied in several works\cite{A,KhS}. It was
essentially solved in \cite{KhS}. The ''mixing properties'' were
found for $g=1$ as a consequence of time delays provided by
saddles. Let us remind that the topology of  typical open
trajectory is the same here as in  the straight line flow (after
removal of domains influenced by the centers). Only saddle points
deform the time functions for the essential part of the
hamiltonian flow for $g=1$. However, the role of these
singularities   is probably completely different for $g=2$. In
this case the topology of foliation was already mixing, so things
like that should essentially remain unchanged.

Now let us consider another interesting example of {\bf The
Maximal Foliations} for $g=4$.

For the maximal  foliations of the class $T$ all boundary saddles
are paired with each other. It means that all of them organize a
system of cycles where the next boundary saddle is paired with the
previous one. The length of every cycle is equal to some even
number $2l_q,q=1,\ldots,f$, so $2l_1+\ldots +2l_f=g$. We say that
the  system has a cycle type $(l_1,\ldots,l_f)$. The maximal
system contains total number of $g$ saddles of the type $<jjkk>$
and $g-2$ saddles of the types $<jklm>$ where all 4 entries are
distinct. For the case $g=4$ we have two possibilities of the
cycle types namely $(1,1)$ and $(2)$. The type $(2)$ is especially
interesting.  This cycle separates a 2-sphere on the South and
North Hemisphere (see Fig 12) where $A_1,A_2,A_3,A_4$ are located
along the equator. The additional two saddles are sitting in the
poles exactly with separatrix curves going to the each ''country''
$A_k,k=1,2,3,4$ along the 4 selected meridians from the north and
south poles.

\begin{lem} The  transitions $A^+_k\rightarrow A^-_l$
with measures $m_{kl}$
 visible from the poles (i.e. located in the corresponding hemisphere)
 are the following:

From the North Pole we can see the transitions $$A_k^+\rightarrow
A^-_l,k=1,3,l=2,4$$
 From the South Pole we can see the transitions
$$A^+_k\rightarrow A^-_l,k=2,4,l=1,3$$ They satisfy to the
Conservation Law $\sum_km_{kl}=|A_l|,\sum_lm_{kl}=|A_k|$, where
$k$ and $l$ are neighbors in the cyclic order $...12341...$. It
implies in particular that $A_1-A_2+A_3-A_4=0$ and ''the constant
flux in one direction of the cycle'' $1234$ is defined provided by
the asymmetry $m$ of the trasition measures $m_{kl}-m_{lk}$ which
is constant:
$$m_{12}-m_{21}=m_{23}-m_{32}=m_{34}-m_{43}=m_{41}-m_{14}=m$$
\end{lem}

\begin{rem}
We say that the  system is rotating clockwise if $m>0$. It is
rotating contrclockwise  if $m<0$.
\end{rem}

\vspace{0.5cm}

\centerline{\Large Section 4. The Homology and Homotopy Classes.}
\centerline{\Large Trajectories and Transversal Curves: $g=2$}

\vspace{0.5cm}

How to describe  the image of transversal semigroups in
fundamental group $\pi_1^+(\Re)\rightarrow \pi_1(V)$ and in
homology group $H_1(V,Z)$?

There are three types of the transversal curves generating all of
them:

(I) The ''Torical Type'' transversal curves not touching the
segments $s$;

(II) The ''Trajectory Type'' transversal curves (or the Poincare
Curves). They coincide with some trajectory of the Hamiltonian
system started and ended in the transversal interval $s^+$ in the
plane $C_1$. We make it closed joining the endpoints by the
shortest transversal interval along $s$;

(III) The general   non-selfintersecting transversal closed
curves.

Let a closed transversal curve $\gamma$ is presented as a point
moving in the combinatorial model. We realize its motion by the
sequence of the almost transversal pieces $\gamma_1\gamma_2\ldots
\gamma_N$. The first piece $\gamma_1\subset C_1$ starts in the
point of the parallelogram $P^1_{0,0}$ just over the segment
$s^+$. It travels in the plane $C_1$ and reaches first time one of
 the segments $s^{\pm}$ located in the parallelogram $P^1_{m_1,n_1}$
of the same plane where $(m_1,n_1)\neq (0,0$. The next path
$\gamma_2\subset C_2$ starts in the corresponding point of the
same segment $s$ but presented  in  another plane $C_2$. We start
it also in the parallelogram $P_{0,0}^2$. It travels in the plane
$C_2$ and ends up in the point of some another segment $s^{\pm}\in
P^2_{m_2,n_2}$ and so on. Finally, the last path $\gamma_N$ ends
up in the same point where the first path started, after  the last
crossing of  some segment $s\subset P^2_{m_N,n_N}$. We may think
that our almost transversal path consists of the pieces of two
kind: the trajectory pieces, passing  streets $p^{\alpha}_k$ from
the initial  segment to the end in positive or negative direction,
and of the orthogonal positively oriented ''jumps'' along some
segment  $s^{\pm}$. Only second type pieces are carrying the
nonzero transversal measure. We can freely create such curves.
Every such curve $\gamma_q$ is almost transversal. However, these
curves are not exactly closed: they start and end up in the
equivalent parallelograms on the equivalent segments $s$ but the
ends are not coincide exactly (even after identification of the
equivalent points). They approach the  same segment $s^{\pm}$ but
(maybe) in the different points. {\bf We can close the ends of the
curve $\gamma_q$ only if the interval from its initial point to
the end one is negatively oriented. } Let us call such almost
transversal curves $\gamma_q$ {\bf semiclosed}.

 Let us present every semiclosed factor $\gamma_q$ with fixed ends
as a product $\gamma_q=\gamma^s_q\tilde{\gamma_q}$ where
$\tilde{\gamma_q}$ is the maximal closed part, and $\gamma^s_q$ is
the ''shortest'' part, i.e. it does not contain any closed piece
inside.  Therefore the curve $\gamma_q$ can be obtained as a
product  $\gamma_q\sim \gamma^s_q\tilde{\gamma_q}$ where
$\tilde{\gamma_q}$ does not cross the segments $s^{\pm}$ at all
(i.e. it is the transversal curve of the Plane Type (I)). Each
path $\tilde{\gamma_q}$ is located completely within the plane
$C_1$ or $C_2$.

Step 1: Classify all closed transversal paths $\tilde{\gamma}\in
C_k,k=1,2$. We demonstrate below the description of their homology
and homotopy types using the ''3-street'' combinatorics described
above.

Let $|a_k|,|b_k|$ denote the measures of the basic transversal
closed curves $a_k,b_k$. {\bf We made our pictures and notations
above such that the cycles $a_k, b_k^{-1}$ are the positive
transversal curves}. According to our construction, the measure
$m$ of the segments $s$ satisfies to the inequality
$$0<m<|a_k|+|b_k|,k=1,2$$ Define the minimal nontrivial nonnegative
pairs of integers $(u_k>0,v_k\geq 0), (w_k\geq 0,y_k> 0)$ such
that
$$m>u_k|a_k|-v_k|b_k|>0$$ and $$m>y_k|b_k|-w_k|a_k|>0$$ (see Fig
13). It means exactly that  these new lattice vectors represent
the shifts $s',s''$ of the segment $s$
 visible directly from $s$ along some shortest trajectories
not crossing other shifted segments, looking to the positive time
direction (+).
\begin{lem}
The  new lattice vectors $h^1_k,h^2_k$ have transversal measures
equal to $u_k|a_k|-v_k|b_k|$ and $y_k|b_k|-w_k|a_k|$
correspondingly. These measures are less than $m$ but their sum is
greater than $m$. Their homology classes
$h^1_k=[a^m_k(+)]=u_k[a_k]+v_k[b_k]$ and $-h_k^2=-[b^m_k(+)]$ for
$[b^m_k(+)]=w_k[a_k]+y_k[b_k]$, generate homological image of the
semigroup of positive closed transversal curves not crossing the
segments $s$ (i.e. of the Torical Type (I)). These classes  have
canonical lifts $a_k^*,b^*_k$ to the free groups generated by
$a_k,b_k$; The lifts $a^*_k,(b^*_k)^{-1}$ of the homology classes
$h_k^1,-h_k^2$ correspondingly generate  semigroup of the homotopy
classes of positive closed transversal curves not crossing the
segments $s_k^{\pm}$, starting and ending in the street number
$0$. These semigroups are free with two generators depending on
the measure $m$ only, whose transversal measures are smaller than
$m$. There exist also the
 similar classes $a^m_k(-),b^m_k(-)$ constructed using the negative
 time direction: their homology classes  are opposite
to the positive ones
$[a^m_k(-)]=-[a^m_k(+)],[b^m_k(-)]=-[b^m_k(+)]$. They represent
the same streets going back. Their lifts $a^m_k(-),b^m_k(-)$ to
the fundamental group $\pi_1(V,s^-_k)$ are defined as a mirror
symmetry of the lifts $a^m_k(+), b^m_k(-)\in\pi_1(V,s^+_k)$ where
$t \rightarrow -t, s^+_k\rightarrow -s^-_k$
\end{lem}
The proof is given below after the reduction to the standard
model.

 Let us make following useful remark:
\begin{lem} Every   path starting and ending in the segment
$s^+_k,k=1,2$ of the same plane $C_k$, has  well-defined homotopy
class in the fundamental group $\pi_1(V)$. Every  curve
 with both ends in any segment of the type $s^{\pm}_k$, has a
well-defined homology class in the group $H_1(V)$.
\end{lem}

Proof. Every segment $s^{\pm}_k$ is realized in the manifold $V$
by the transversal segment joining two saddle points. After
cutting along these segments, we split our surface into two
disjoint pieces. Every piece has a boundary
$$a^*_1b^*_1(a^*_1)^{-1}(b^*_1)^{-1}\sim \kappa=s^+_k\bigcup
s^-_l,k\neq l$$. The transversal segments are identified with each
other according to the rule $$s^+_1=s^-_2,s^+_2=s^-_1$$ So the
homology 1-classes modulo boundary $H_1(V,\kappa)$ are the same as
in $H_1(V)$. We may assign  homology class to every piece of
trajectory passing any street from the beginning to the end.

Concerning homotopy classes, we assign the invariant
$\phi_{\alpha\beta}\in\pi_1(V,s_1^+)$ to every piece of trajectory
passing two streets $$\gamma_q\subset
p^{\beta}_2p^{\alpha}_1=\{\alpha\beta\},\phi:\gamma_q\rightarrow
\pi_1(V,s^+_1)$$ We use here the fact that the transversal segment
$s^{\pm}$ is contractible. Therefore  the choice of initial point
in it is unimportant. Every infinite trajectory $\gamma$ is coded
by the infinite sequence of pieces $$\gamma=\ldots
\{\alpha_q\beta_q\}\{\alpha_{q-1}\beta_{q-1}\}\ldots$$ So the can
apply the homomorphism $\phi$ into the fundamental group for all
sequence. Every finite connected even piece $\delta'$ of this
sequence define the element $\phi(\delta')\in \pi_1(V,s^+_1)$. Its
ends can be joined by the shortest transversal piece along the
segment $s$. Depending on orientation of this piece, either
$\phi(\delta')$ or $\phi(\delta')^{-1}\in\pi_1(V,s^+_1)$ define a
closed positive transversal curve.

{\bf Reduction to the Standard Model}:  The segment $s$ of the
length $m$ is divided on  $5$ connected open pieces $\tau_q$:
There exist
 exactly 9 possible types of trajectory pieces
$\{\alpha\beta\},\alpha,\beta=1,0,2$ with measures
$p_{\alpha\beta}$ but 4 of them are in fact empty. They   have
 measure equal to zero. In order to see that, we remind how these pieces
 were constructed. A segment $s=s_1$ of the total measure $m$ is
 divided into 3 pieces by the points $0,1,2,3$ for $k=1$ and by
 the points $0',1',2',3'$ for $k=2, s=s_2$: For the streets  we have:
 $$ \alpha=1,0,2=[01],[12],[23],\beta=1',0',2'=[0'1'],[1'2'],[2'3']$$
 So,  the index $\alpha=1,0,2$ corresponds
 to the segments $[01],[12],[23] $ with measures $p^{\alpha}_1$,
 and the index $\beta=1',0',2'$  corresponds
 to the segments $[0'1'],[1'2'],[2'3']$ with measures $p^{\beta}_2$.
The positions of the points $0=0',3=3'$ are fixed. Other points
never coincide for the generic foliations.

 Every jump from $C_1$ to $C_2$  is accompanied by the permutation
of 3 segments: the left street number 1 in $C_1$ ends up in the
extreme right part of $s=s'$ before making jump to $C_2$. The
right street number 2 in $C_1$ ends up in the extreme left part of
$s=s''$ in $C_2$. So jumping from $C_1$ to $C_2$, we should
permute the segments $1=[01]=1$ and $2=[23]=2$ preserving
orientation:
$$\eta_{12}:2\rightarrow 0=2^*,
 3\rightarrow 3^*\in s,
 1\rightarrow 3=1^*, 0\rightarrow 0^*\in s$$.
 Here $ |[23]|=|[2^*3^*]|,|[01]|=|[0^*1^*]|$.
The segment $s$ is divided by the points
 $0=2^*,1',2',3^*,0^*,3=1^*$, so it is presented as a    union of   5 sub-segments
$$s=\tau_1+...+\tau_5$$
 In order to return
back from the second torus (plane) $C_2$ back to $C_1$, we need to
apply the similar map $\eta_{21}$ based on the permutation of the
streets $2'=[2'3']$ and $1'=[0'1']$. Our broken isometry
$i_{\sigma}$ based on the permutation $\sigma$  of 5 pieces, is
defined as a composition
$$i_{\sigma}=\eta_{21}\eta_{12}$$
 \begin{lem}6 possibilities called {\bf the Topological Types of Foliations},
 for the nonzero set of 5 measures
 $p_{\alpha\beta}$ exist here (the measures $p_{\alpha\beta}$ of the
sub-segments $\tau_q,q=1,2,3,4,5$ are given
 in the natural order on the segment $s$):
$$(I):0=2^*<1'<2'<3^*<0^*<1^*=3 ; \sigma=(32541);$$
$$p_{12}+p_{02}+p_{21}+p_{20}+p_{22}=m$$
$$(II):0=2^*<1'<3^*<2'<0^*<1^*=3; \sigma=(24153);$$
$$p_{12}+p_{01}+p_{02}+p_{21}+p_{20}=m$$
$$(III):0=2^*<1'<3^*<0^*<2'<1^*=3;\sigma=(41523);$$
$$p_{10}+p_{12}+p_{00}+p_{21}+p_{20}=m$$
$$(IV):0=2^*<3^*<1'<2'<0^*<1^*=3;\sigma=(25314);$$
$$p_{12}+p_{01}+p_{00}+p_{02}+p_{21}=m$$
$$(V):0=2^*<3^*<1'<0^*<2'<1^*=q3;\sigma=(31524);$$
$$p_{10}+p_{21}+p_{01}+p_{00}+p_{21}=m$$
$$(VI):0=2^*<3^*<0^*<1'<2'<1^*=3;\sigma=(52134);$$
$$p_{11}+p_{10}+p_{12}+p_{01}+p_{02}=m$$ All other measures
$p_{\alpha\beta}$ are equal to zero. We have
$$\sum_{\alpha}p_{\alpha\beta}=p_2^{\beta};
\sum_{\beta}p_{\alpha\beta}=p_1^{\alpha}$$
\end{lem}

 As we can see, our 3-street model automatically
creates the standard type broken isometry
$\eta_{21}\eta_{12}=i_{\sigma}:s\rightarrow s$ generated by the
permutation $\sigma$ of 5 pieces $\tau_q$ of the segment $s$.
 Such systems were studied by the ergodic people since 1970s.  Let
us point out that our combinatorial model easily provides  full
information about the topology lying behind this permutation. We
know geometry of all pieces. Let us define also a {\bf Shift
Function}: The  map $i_{\sigma}:\tau_q\rightarrow  s$ for every
sub-segment $q=1,2,3,4,5,$ is an orientation preserving isometry
inside. So it is a shift $x\rightarrow x+ r_q$ inside of these
strips.

 We call   $r_j\in R$ a value of the {\bf  Shift Function}. Iterating our
system, we have  similar ''Shift Function'' for every piece of the
trajectory $\gamma$ if it started and ended in the transversal
segment $s^+$. We call $\gamma$ {\bf a Positive Piece} if its
Shift Function $r(\gamma)$ is positive. In the opposite case it is
{\bf a Negative Piece}. For the irreducible generic systems there
are no periodic solutions (i.e. no pieces with zero shift
functions). Every piece defines {\bf a Closed Transversal Curve of
the Trajectory Type} closing it by the shortest path along $s$.
Negative pieces define positive closed transversal curves and vice
versa. Let us calculate  the value of the shift function:

Let us define an algebraic object (no doubt,  considered by the
ergodic people many years ago): {\bf an Associative Semigroup}
$S_{\sigma,\tau}$ with ''measure''. It is generated by the 5
generators $R_1,...,R_5\in S_{\sigma,\tau}$. There is also a zero
element $0\in S_{\sigma,\tau}$. We define multiplication in the
semigroup $S_{\sigma,\tau}$ as in the free one but some ''zero
measure''  words are equal to zero by definition. In order to
define which words are equal to zero, we assign to every generator
an interval $R_q\rightarrow \tau_q,q=1,...,5$.

{\bf  The word $R=R_pR_q$ is equal to zero in the semigroup
$S_{\sigma,\tau}$ if and only if $\tau_p\bigcap i_{\sigma}(
\tau_q)=\emptyset$.} We assign to the word $R=R_pR_q$ the set
$\tau_R$ such that
$$i_{\sigma}(\tau_{R_pR_q})=i_{\sigma}(\tau_R)=\tau_p\bigcap
i_{\sigma}(\tau_q)$$ if it is nonempty. By induction, we assign to
the word $R_pR$ the set $\tau_{R_pR}$ where
$$\tau_{R_pR}=i_{\sigma}^{-N}(\tau_p)\bigcap \tau_R$$ if it is
nonempty. We put $R_pR=0$ otherwise. Here $N$ is the ''length'' of
the word $R$. This semigroup is associative. For every word in the
free semigroup $R=R_{q_1}...R_{q_N}$ we assigned a set
$$i_{\sigma}^{-N+1} (\tau_{q_1})\bigcap
i_{\sigma}^{-N+2}(\tau_{q_2})\bigcap...\bigcap \tau_{q_N}=\tau_R$$
whose measure is well defined. We put $R$ equal to zero if the
measure is equal to zero.

The ordered  sequence $R_{\infty}=\prod_{p\in Z} R_{q_p}$,
infinite in both directions, defines trajectory of the flow if and
only if  every finite sub-word is nonzero. The measure is defined
for the ''cylindrical'' sets $U_R$ consisting of all
''trajectories'' with the same sub-word $R$ sitting in the same
place for all of them. It is equal to the ''measure'' of the word
$\tau_R$. The ''Shift Function'' is also well-defined by the
semigroup and every trajectory. It does not define a homomorphism
because it can be nonzero for the word equivalent to zero. For the
nonzero words its values are bounded $r(R)<m, 0\neq R\in
S_{\sigma,\tau}$.

{\bf We present below the calculation of representation of this
semigroup in the fundamental group of the Riemann Surface
generated by the real part of the holomorphic one-forms.}

 If we know  the time characterization functions in
each piece $\{\alpha\beta\}$ (it is simply  a sum of times in
every street $p^{\alpha}_1,p^{\beta}_2$), we can also study the
Hamiltonian system with corresponding natural time. All
trajectories $\gamma\in\{\alpha\beta\}$ have the same homotopy
classes $\phi_{\alpha\beta}\in\pi_1(V)$ known to us (see below).
 The global reduction of the flow
is based on the non-closed transversal Poincare section, i.e. on
the transversal segment $s$ leading from one saddle to another.
This construction seems to be  the best possible genus 2 analog
 of the reduction of straight line flow on the 2-torus  to the
rotation of circle.  Many features of this construction certainly
appeared
 in some very specific examples
 studied before.

We are going  to classify now the homotopy and
homology classes of the trajectories  starting and
ending in the transversal segments $s$. Before doing that, we need
to proof  the lemma above in order to finish the first step..

Proof of the lemma. We denote the streets $p^{\alpha}_k$ in the
plane $C_k$ simply by the symbols $\alpha=0,1,2$ where the longest
one is $0$. It is located between two others. The right one in
$2$, and the left one is $1$. For the sum of their homological
''lengths'' in the simplified notations $h^{\alpha}_k\rightarrow
[\alpha]$ we have $[1]+[2]=[0]$. This relation will be treated
also homotopically later. For the widths (i.e. transversal
measures) we have $\sum_{\alpha}|p^{\alpha}_k|=m$. Their bottom
parts cover together the segment $s$. This picture is invariant
under the change of direction of time and simultaneous permutation
of the lower and upper segments $s$. All transversal paths  in the
plane $C_k$ can be written combinatorially in the form $$\ldots
\rightarrow\alpha_q \rightarrow 0\rightarrow
\alpha_{q+1}\rightarrow 0\rightarrow\alpha_{q+2}\rightarrow
0\rightarrow \ldots$$ where $\alpha_m=1,2$. We have a pair of
positive ''basic cycles'' $[0\rightarrow 1\rightarrow
0]=(b^m_k(+))^{-1}$ and $[0\rightarrow 2\rightarrow
0]=a^m_k(+),k=1,2$. All other transversal cycles not touching the
segment $s$ and its shifts, starting and ending in the longest
strip $0$, have a form of the arbitrary word in the free semigroup
generated by $a^*_k=a^m_k(+),(b^*_k)^{-1}=(b^m_k(+))^{-1}$ in the
plane $C_k$. The measures of these new $m$-dependent basic cycles
are
$$|a^*_k|=|a^m_k|=|p^0_k|+|p^2_k|,|b^*_k|=
|b^m_k|=|p^0_k|+|p^1_k|,k=1,2$$. Topologically the cycles
$a^m_k(+),b^m_k(+)$ represent some canonical
 $m$-dependent basis of
 elementary shifts in the group $Z^2_k$. We choose the simplest transversal
paths joining the initial point with its image  in order to lift
this basis  to the free group. These shifts map the segment
$s\subset P^k_{0,0}$ exactly into the the segments $s',s''$
attached to the middle part of the long street from the left side
(for $a^m_k(+)$) and from the right side (for $b^m_k(+)$)--see Fig
13.
 The subgroup generated by $a_k,b_k$ in fundamental group is free.
  Homologically  the
elements $[a^*_k]=h^1_k=a^m_k(+),[b^*_k]=h^2_k=b^m_k(+)$
 are calculated in
the formulation of lemma using the transversal measures
$|a_k|,|b_k|$ of basic cycles and the measure $m$. From geometric
description of new cycles $a^*_k,b^*_k$ as of the paths we can see
that they satisfy to the same relation as the original $a_k,b_k$:
Their commutator path exactly surrounds one segment $s$ on the
plane. Lemma is proved. This is the end of the Step 1.

We always choose initial point on the segment $s^+_1$ in the plane
$C_1$. It is the same as to choose initial point on the segment
$s^-_2=s^+_1$ in
 the plane $C_2$.  We are going to use a  new basis
 of fundamental group $\pi_1(V,s^+_1=s^-_2)$:  {\bf The New $m$-Dependent
Transversal Canonical Basis} (see Fig 13) is
 $$a^*_1=a^m_1(+),b^*_1=b^m_1(+),a^*_2=a^m_2(-),b^*_2=b^m_2(-)$$
 We attach this basis  to
 the segment $s^+_1=s^-_2$ in order to
treat them as the elements of fundamental group $\pi_1(V,s^+_1)$.
\begin{lem} The homology classes $h^{\alpha}_k$ of the  streets
number $\alpha$ in the plane $C_k$,  are following:
$$h^1_k=[a^m_k(+)],h^2_k=[b^m_k(+)],h^k_0=[a^m_k(+)]+[b^m_k(+)]\in
H_1(V,Z)= H_1(V,s^+\bigcup s^-)$$ Here $k=1,2$. We assume that the
streets are oriented to the positive time direction. For the
negative time we simply change sign
$[a^m_k(-)]=-h^k_1,[b^m_k(-)]=-h^k_2$.
\end{lem}
This lemma immediately follows from the description of the
homology classes $[a^m_k(\pm)]$ and $[b^m_k(\pm)]$.

 Let us describe the homotopy classes
$\phi_{\alpha\beta}\in \pi_1(V,s^+_1)$
 of the two-street paths $\alpha\beta$ starting and ending on the same open segment
$s=s^+_1=s^-_2 \in C_1\bigcap C_2$. Here $C_1\bigcap C_2$ is a
union of all shifts of the cycles $s^+\bigcup s^-$.

 The streets  in  the plane $C_k$ end up in the points of the segment $s^-_k$.
 Consider first the plane $C_1$.
 In order to make closed paths out of the streets $p^{\alpha}_1$,
 we need to extend them: we go
 around the segment $s$ from $s^-_1$ to $s^+_1$ from the
 right or from the left side of it (see Fig 14):

For the representation of  street number 1 we use  the  path
$a^*_1=a^m_1(+)$ closed by passing $s'$ from
 the right side along the path $\kappa_1$, circling contr-clockwise around $s'$;
For the street number 2 we use  the path $b^*_1=b^m_1(+)$ closed
by passing $s''$ from
 the left side along the path $\kappa_2$ circling clockwise around $s''$;
 So we have $$p^1_1\sim a^*_1\kappa_1^{-1}$$
$$p^2_1\sim b^*_1\kappa_2^{-1}$$ Here and below {\bf the symbol
$\sim$ means ''homotopic with fixed ends''}.
 We define a closed path circling clockwise around the segment $s$ in $C_1$:
 $$\kappa=\kappa_2\kappa_1^{-1}\sim (s^+\bigcup s^-) \sim
a^*_1b^*_1(a^*_1)^{-1}(b^*_1)^{-1}$$ For the street number 0 we
assign  the path $$p^0_1\sim a^*_1b^*_1\kappa_2^{-1}$$ The same
description we have also for the streets
$p^{\alpha}_2(-)=(p^{\alpha}_2)^{-1}$ in the plane $C_2$ going
 in the opposite direction, replacing (+) by (-) and
the paths $\kappa,\kappa_1,\kappa_2$ by the similar paths $\delta,
\delta_1,\delta_2$: $$(p^1_2)^{-1}\sim
a^*_2\delta_1^{-1},(p^2_2)^{-1}b^*_2\sim\delta_2^{-1},
(p^0_2)^{-1}\sim a^*_2b^*_2\delta_2^{-1}$$. We have
$$\kappa_1=\delta_2,\kappa_2=\delta_1,\delta=\delta_1\delta_2^{-1}=\kappa^{-1}$$
We use the new basis $a^*_1,b^*_1,a^*_2,b^*_2$ defined above,
dropping the measure $m$ and signs $\pm$, as it was indicated
above.  {\bf All Formulas are written in the new basis}.

We are ready to present  the  two-street formulas in $\pi_1$,
where at first we are passing the street $\alpha$ in $C_1$ and
after that the street $\beta \in C_2$. Writing $\alpha\beta$
instead of $\phi_{\alpha\beta}$ and performing very simple
multiplication of paths, we obtain following formulas in the group
$\pi_1(V,s^+_1)$:

\begin{lem}
The homotopy types of all nonnegative almost transversal
two-street passes in the positive time direction starting and
ending in the segment $s^+_1=s^-_2$,  including the trajectory
passes, are equal to the following list of values of the elements
$\phi_{\alpha\beta}\in \pi_1(V,s^+_1)$ where $\alpha=1,0,2$ and
$\beta=1',0',2'$:

 $$11'\sim a^*_1\kappa^{-1}(a^*_2)^{-1}; 10'\sim a^*_1b^*_1(a^*_2)^{-1};
  12'\sim b^*_1(a^*_2)^{-1}$$

$$01'\sim a^*_1(b^*_2)^{-1}(a^*_2)^{-1}; 00'\sim a^*_1b^*_1\kappa
 (b^*_2)^{-1}(a^*_2)^{-1}; 02'\sim b^*_1\kappa(b^*_2)^{-1}(a^*_2)^{-1}$$

$$21'\sim a^*_1(b^*_2)^{-1}; 20'\sim a^*_1b^*_1\kappa(b^*_2)^{-1}; 22'\sim
b^*_1\kappa(b^*_2)^{-1}$$
 For the negative time direction the two-street passes have homotopy classes equal to
$\phi_{\alpha\beta}^{-1}$ in the same group
$\pi_1(V,s^-_2)=\pi_1(V,s^+_1)$

 For the
Topological Types I--VI of Foliations following homotopy classes
of  almost transversal  two-street passes have nonzero measure:
$$\phi_{\alpha\beta}=p^{\alpha}_1xp^{\beta}_2,\phi^*_{\alpha\beta}=
(p^{\alpha}_1)^{-1}x(p^{\beta}_2)^{-1}$$  (here $x\geq 0$ means
positive transversal shift along the segment $s$):

 Type (I): All 9 classes $\phi_{\alpha\beta}$, and all  classes
 $\phi^*_{\alpha\beta}$ except $11',10',01',00'$

Type (II): All classes $\phi_{\alpha\beta}$ except $22'$, and all
classes $\phi^*_{\alpha\beta}$ except $11',10',01'$

Type (III): All classes $\phi_{\alpha\beta}$ except $22',02'$, and
all classes $\phi^*_{\alpha\beta}$ except $11',01'$

Type (IV): All classes $\phi_{\alpha\beta}$ except  $20',22'$, and
all classes $\phi^*_{\ alpha\beta}$ except $11',10'$

Type  (V): All classes $\phi_{\alpha\beta}$ except $22',20',02'$,
and all classes $\phi^*_{\alpha\beta}$ except $11'$

Type (VI): All classes $\phi_{\alpha\beta}$ except
$20',22',00',02'$, and all 9 classes $\phi^*_{\alpha\beta}$

\end{lem}

Next theorem follows  the standard scheme of the ergodic theory,
but the values of topological quantities are explicitly calculated
within the combinatorial model of the flow on the Riemann Surface:
\begin{thm} The combinatorial model of foliation defines an Ensemble
 consisting of the following ingredients:

1.The semigroup $S_{\sigma,\tau}$ and its representation in the
fundamental group are given. All
 positive finite words $$R=R_{j_1}R_{j_2}\ldots R_{j_N}\in \pi_1(V)$$
are written in the new $m$-dependent  Transversal Canonical Basis
$a^*_1,b^*_1,a^*_2,b^*_2$, of all
 lengths $||R||=N\geq 1$,
where every symbol $R_j$ is equal to one of the elements
$\phi_{\alpha\beta}\in \pi_1(V),j=1,2,3,4,5$. 2.The  permutation
$\sigma$ and measures of the sub-segments are given. We have
$\tau_q,q=1,2,3,4,5, \sum_q\tau_q=m$. They  represent the pairs
$\alpha\beta$ with nonzero measure,  according to the types
$I,II,III,IV,V,VI$ above. The permutation $\sigma$ defines a
broken isometry $i_{\sigma}:s\rightarrow s$ well-defined for the
inner points of the sub-segments $\tau_q$.

2.For every length $N$, the set of words with nonzero measure is
ordered in the following way: Every such word $R\in \pi_1(V)$ is
represented by a single connected sub-segment $tau_R\subset s$
with nonnegative length. Every word represented by the empty
sub-segment $\tau_R=\emptyset$, is equivalent to zero.  Nonzero
segments do not intersect each other for $R\neq R'$ . They are
naturally ordered in the segment $s$ of the length $m$, and sum of
their measures is eq1ual to $m$ for every length $N$. In order to
multiply any word $R$ from the left by the elementary word
$R_q=\phi_{\alpha\beta},q=1,2,3,4,5$, we apply the map
$i_{\sigma}^{-N+1}$ to the set $\tau_q$ segments of the nonzero
measure and intersect it with segment $\tau_R$ representing $R_j$:
$$\tau_{R_qR}=i_{\sigma}^{-N+1}(\tau_{q})\bigcap \tau_R$$ If
intersection is empty, the product $R_qR$ is equivalent to zero by
definition. The set of words $S$ factorized by the words
equivalent to zero, forms an associative semigroup
$S_{\sigma,\tau}$ with $0\in S_{\sigma,\tau}$ defined by our
foliation $dH=0$. If intersection
$i_{\sigma}^{-N+1}(\tau_q)\bigcap \tau_R$ is nonempty, we assign
it to the corresponding product word $R_qR$ of the length $N+1$. A
width of this strip is treated as a nonzero measure of the new
word. The Shift Function $r(R)$ for every nonzero word
$R=R_{j_1}\ldots R_{j_N}$ is equal to the sum $$r(R)=\sum_q
r(R_{j_q})$$ It is bounded for all nonzero words $r(R)<m$.

3.Every finite word  $R\in S_{\sigma,\tau}$ defines a connected
strip of the nonseparatrix  trajectories of the Hamiltonian System
$dH=0$ corresponding to our data, with the transversal measure
$\tau_R$, starting and ending in the transversal segment $s^+_1$.
It defines a positive closed transversal curve $\gamma_R$ with
transversal measure equal to $r(R)$ if $r(R)<0$, and negative
closed transversal curve if $r(R)>0$. The transversal measure of
these curves are equal to $-r(R)$. The individual infinite
trajectories are presented by the infinite sequences of the
symbols $R_j$ such that every finite piece of the sequence is
nonzero as an element of the semigroup $S_{\sigma,\tau}$. The
measure on the set of trajectories is defined by the transversal
measure of the strips corresponding to the finite words in $S$.
The basic measurable sets are ''cylindrical'' (i.e. they consists
of all trajectories with the same finite word $R$, and measure is
equal  to $\tau_R$).

\end{thm}

{\bf Problem 1}: How to calculate effectively the  words in the
free group with 2 generators $a,b$ describing new $m$-dependent
canonical transversal basis $a^*,b^*$ in the 2-torus with
straight-line flow and obstacle $s\subset T^2$ of the transversal
size $m$? We assume that $m\rightarrow 0$, and the initial
transversal canonical basis $a,b$ remains fixed.

Therefore the Step 2 is realized. We described  the homology and
homotopy classes associated with trajectories starting and ending
in the segments $s$. They   correspond to the ''Trajectory Type ''
closed transversal curves.

Now let us start to discuss the most complicated

{\bf Problem 2}: How to describe all transversal curves?

We are going to study this problem extending the method of the
previous section where the trajectory type curves were described.

It is easy to prove following statement:

\begin{lem} Consider the transversal homotopy classes (with fixed ends
$\gamma(0)$ and $\gamma(1)$) of positive almost transversal curves
$\{\gamma\}$ starting and ending in the segment $s=s^+_1=s^-_2$.
Every such curve within its class can be presented by the almost
transversal path $\gamma$ consisting of the  following pieces:
1.Any full street $(p^{\alpha}_k)^{\pm 1}\in C_k,k=1,2$ passed in
any direction;  2.Any jump $x\in (0,m)$ along the segments
$s^{\pm}_k,k=1,2$ in positive direction before passing the next
street. Therefore every such curve can be presented by the word
consisting of these symbols, starting and ending by the streets:
$$W'= (p^{\alpha_1}_{k_1})^{\pm
1}x_{1}(p^{\alpha_{2}}_{k_{2}})^{\pm 1}\ldots
x_{N-1}(p^{\alpha_{N}}_{k_{N}})^{\pm 1}$$ The number $N$ is even.
These pieces are satisfying to the following obvious relation: If
the same street was passed in some direction and back immediately
after that (maybe with some positive jump $x$ between them), it
can be removed, and the jump between them can be divided into two
parts $x=x'+x''$ added to the previous and to the next steps.

\end{lem}

By the {\bf Topological Type $W$ of any Word $W'$} we call the
same word where all transversal jumps $x_q$ are simply omitted. We
say that $(x_1,\ldots, x_{N-1})$ belongs to the {\bf Existence
Domain of the given Topological Type} $\Delta_W$ if for all jumps
belonging to this domain there exists an almost transversal curve
of that type. We assign to every such Topological Type a positive
measure equal to the measure of the domain $\Delta_W$ in the cube
$[0m]^{N-1}$. This measure depends on the ends of the curve. We
can close this  curve and construct a closed transversal positive
curve only if the last end $\gamma(1)$ is located from the left
side of the initial point $\gamma(0)$. Taking
$x_N=|[\gamma(N),\gamma(0)]|$. The transversal measure of this
closed transversal curve is equal to
$\oint_{\gamma}dH=\sum_{q=1}^Nx_q$.

We assign  homology class to every transversal curve $\gamma$ with
ends in the same segment $s^+_1=s^-_2$: $$[\gamma]=\sum
(\pm)[p^{\alpha_q}_{k_q}]\in H_1(V,Z)$$ depending on the
topological type only. The homology classes corresponding to the
streets were calculated above.

In order to calculate the homotopy classes, we consider two-street
passes as above. The list of  homotopy classes of two-street
passes in the group $\pi_1(V,s^+_1)$ was done before, for the case
where both of them are positively directed (or both are negatively
directed like trajectories).For the  trajectories, without
transversal jumps,  we had only 5 types of the trajectory type
two-street passes $\phi_{\alpha\beta}$ with nonzero measures. Here
we have more as it was indicated in the lemma above. Besides that,
we may have two-street transversal passes concentrated in one
plane, of the types
$$\psi^{\alpha\alpha'}_1=p^{\alpha}_1x(p^{\alpha'}_1)^{-1},x\geq 0 $$
$$\psi^{\beta\beta'}_2=(p^{\beta}_2)^{-1}xp^{\beta'}_2, x\geq
0$$ These transversal cycles in fact were described above: they
form the new $m$-dependent Transversal Canonical Basis.

\begin{lem}

Following positive transversal two-street passes of the types
$\psi^*{\alpha\alpha'}_1$, starting and ending in the segment
$s^+_1=s^-_2$, have nonzero measure only:
$$\psi^{01}_1=p^0_1x(p^2_1)^{-1}\sim a^*_1;
\psi^{10}_1=p^1_1x(p_0^1)^{-1}\sim (b^*_1)^{-1}$$
$$\psi^{12}=p^1_1x(p^2_1)^{-1}\sim (b^*_1)^{-1}a^*_1$$ For the
two-street passes of the type $\psi^{\beta\beta'}$  we have
$$\phi^{0'1'}=a^*_2,\phi^{2'0'}=(b^*_2)^{-1};
\phi^{2'1'}=(b^*_2)^{-1}a^*_2$$ Here $1,0,2$ are the numbers of
streets for $k=1$ and   $1',0',2'$ are the numbers of streets for
$k=2$. Combining these formulas with homotopy classes of positive
and negative almost transversal two-street passes
$\phi_{\alpha\beta}$ and $\phi^*_{\alpha\beta}$ calculated above,
we have a complete list of two-street transversal passes of all
types jumping in positive transversal direction.
\end{lem}
For every transversal curve $\gamma$ we consider now the
deformations within its topological type such that we allow to
move the first end $\gamma(0)$  to the left along the segment $s$
(the last end $\gamma(1)$ is staying). We allow to move left all
streets passes except the last one  changing the jumps $x_q$.
Finally we are coming to the following.  We find an extremal
(minimal) transversal distance between $\gamma(0)$ and
$\gamma(1)$. If there is a
 position (within the given topological type) such that
 the segment $\tau$ starting in $\gamma(1)$ and
ending in $\gamma(0)$, is positive, we say that this topological
type represents a positive closed transversal curve
$\bar{\gamma}=\gamma\bigcup \tau$.

 {\bf  Every positive closed transversal curve $\gamma$ starting
 and ending in the segment $s^+_1=s^-_2$, has a homotopy class
  which is a positive product
of the following elements in the group $\pi_1(V,s^+_1)$}:
$$\phi_{\alpha\beta},\phi^*_{\alpha\beta},\psi^{\alpha\alpha'},
\psi^{\beta\beta'}$$ However, not all of them represent any closed
transversal curve:

1.The measure of this type should be nonzero. 2. The initial point
$\gamma(0)$ for some representative should be located to the left
from the endpoint $\gamma(1)$. The algorithm for finding such
representative (or testing whether it exists or not)  consists of
the ''left'' deformations of street passes changing the jumps
$x_q$, starting from the initial point $\gamma(0)$, and then
repeating this process back starting from the motion of endpoint
$\gamma(1)$ to the right.

\vspace{0.3cm}

{\bf How to classify all non-selfintersecting transversal closed
curves?}

Let us describe first Solution of this Problem for the curves not
touching the segment $s$, i.e. belonging to one plane $C_k$.

Consider first a 2-torus without one point $T^*=T^2 (minus) P$
presented as a standard parallelogram $P\subset C$ with one inner
point missing, $\partial P=aba^{-1}b^{-1}$.

following facts are well-known:

1.Every non-selfintersecting closed curve $\gamma\subset T^*$
representing  homology class $[\gamma]=ka+l(-b)\subset
H_1(T^2,Z)=H_1(T^*,Z)$ where $(k,l)=1$, can be deformed in $T^*$
to the  standard form such that it has $k$ transversal
intersection points with the cycle $b$, and $l$ transversal
intersection points with $a$.

 Let us choose our notations for the basic cycles $a,b$
 and direction of curve
$\gamma$ such that $k>|l|>0$.  There are two cases here: (I)$l>0$
and (II)$l<0$. There are also trivial cases such that either
$k=0$, or $l=0$, or $|k|=|l|=1$.

 2.Following algorithm allows to write effectively the presentation
 of a positive element $k>0.l>0$ as a positive word in the alphabet
  $a=a,b=b^{-1}$, unique up to  the cyclic permutation:

 {\it Let our curve $\gamma$ starts at the vertex of the
parallelogram $P$. It is represented in $P$ by the sequence of
$k+l$ segments $t_j\subset P,j=1,...,k+l-1$, starting and ending
in the crossing points of the curve $\gamma$ with the boundary of
parallelogram $P$: Consider only the case (I). We choose notations
for these points such that $y_1,y_2,\ldots, y_l\subset a$ and
$y_{l+1},\ldots, y_{k+l}\subset b$ and equivalent points $y'_j$ in
the components $a^{-1},b^{-1}$. Here the points $y_1$ and
$y_{k+l}$ are the opposite vertices of $P$ along the diagonal, all
other points are distinct. All segments  have positive length. We
simply denote this crossing points by their   numbers in the lemma
below $y_j\rightarrow j,y'_j\rightarrow j'$ (see Fig 15). The
topological type of  non-selfintersecting closed curve $\gamma$ in
$T^*$ is completely determined by its homology class and by the
number of the domain $S_j$ between two neighboring segments
containing the removed point $P\in T^2 $. The list of segments
beginning from the left upper vertex $a\bigcap b\in P$, naturally
ordered by the position in the Parallelogram $P$, is following for
the Case (I): $$t_1=[l,l+1],t_2=[l-1,l+2],..., t_l=[1,2l]$$
$$t_{l+1}= [(l+1)',2l+1],...,t_{k}=[k',(k+l)]$$
 $$
t_{k+1}=[(k+1)',l'],..., t_{k+l-1}=[(k+l-1)',2']$$ These segments
divide the open parallelogram $P$ into $k+l$ domains $$P
minus(\bigcup_j t_j)=S_1\bigcup S_2\ldots \bigcup S_{k+l}$$ where
$S_j$ is bounded by  $t_j$, $t_{j-1}$ and $\partial P$.

Every segment $t_j$ divides the parallelogram $P$ into the upper
and lower parts  $P minus (t_j)=P_j^+\bigcup P^-_j$ where
$t_q\subset P^+_j$ for $q<j$. Fixing the domain $S_r$ containing
the removed point, we see that  all domains $P^+_j$ for $j\geq r$
also contain this point. We call them the {\bf Marked Domains}.
All domains $P^+_j$ with $j<r$ are not marked, i.e. they do not
contain this point.

Consider now a naturally ordered sequence of the domains
$P^+_{q_1}P^+_{q_2}...P^+_{q_{k+l-1}}$ along the closed
non-selfintersecting curve $\gamma\subset C$ not crossing the
removed point and its shifts. The pieces $t_j\subset \partial
\bar{P}^+_j$ exactly form a closed curve $\gamma=\bigcup_j t_j$.
We start from the point $y_1$, i.e. from the domain
$P^+_l=P^+_{q_1}$, and end up with the domain $P_k=P_{q_{k+l}}$
ending in the point $y_{k+l}$.

Assign now to every domain $P^+_j$ following element of the free
group with two generators:

$$\psi:P^+_j\rightarrow \kappa^c; \kappa=b^{-1}aba^{-1}$$ if the
boundary of its closure $\partial\bar{P}^+_j$ does not contain
anyone of the full cycles $a$ and $b$

$$\psi:P^+_j\rightarrow a\kappa^c $$ if the boundary of its
closure contains only a full cycle $a$, not $b$

$$\psi:P^+_j\rightarrow \kappa^c b$$ if the boundary of its closure
contains only a full cycle $b$ but not $a$

$$\psi: P^+_j\rightarrow a\kappa^c b$$ if the boundary of its
closure contains fully both cycles $a,b$

Here $c=0$ if our domain $P^+_j$ is non-marked, and $c=1$ if it is
marked.

The products of elements of free group with two generators
$$\psi(\gamma)=\psi(P^+_{q_1})...\psi(P^+_{q_{k+l-1}})$$ is equal
to the homotopy class of the closed curve $\gamma$ in the free
group with two generators. We call this product {\bf An Upper
Triangle Decomposition} of the class $[\gamma]\in \pi_1(T^*)$ This
description gives a complete classification of the homotopy
classes of non-selfintersecting closed curves in the punctured
2-torus $T^*$, representing the homology class $$[\gamma]=ka+lb\in
H_1(T^*,Z)=H_1(T^2,Z), k>l>0, (k,l)=1$$ The homology class and the
integer number $1\leq r\leq k+l$ completely determine the homotopy
class of the nonselfintersecting closed curve in $T^*$.

Therefore every indivisible positive homology class $k[a]+l[-b]\in
H_1(T^*,Z),k>0,l>0$, can be represented by one positive word of
the length $k+l$ in the free group with two positive generators
$a=a',b^{-1}=b'$ realized by the non-selfintersecting curve. This
word is unique up to cyclic permutation (a source for the number
$r$).}

Let us assume now that we have a straight line flow on the 2-torus
with irrational rotation number (i.e. it is generic) and fixed
transversal canonical basis $a,b$ where $a,b^{-1}$ are positive as
above. Every indivisible homology class in $H_1(T,Z)$ can be
realized by the nonselfintersecting closed curve transversal to
foliation. We can see that {\bf semigroup of positive closed
transversal curves has infinite number of generators containing
the elements with arbitrary small transversal measure}. The same
result remains true after removal of one point $T\rightarrow T^*$
but the semigroup became nonabelian. Remove now a transversal
segment $s$ with positive measure $m>0$ from the torus
$T^*\rightarrow T^{(m)}$. As we shall see below, the situation
changes drastically: the semigroup became finitely generated (in
fact, there are two generators in it).

Consider now the {\bf Classification Problem} of the transversal
closed curves for the hamiltonian foliations with Transversal
Canonical Basis which do not touch  the segment $s$. They are
concentrated in the flat part $T^*_m=T^2 minus (s)$  in one plane
$C_k$ only. Let $k=1$ (we drop the number $k$ beginning from now).
It is convenient to use a three-street model for this goal. As
above, we denote the maximal street by the number $0$, the left
one by $1$ and the right one by $2$. Our Fundamental Domain
consists of the union of these three streets. We use an extension
of it adding 2 more streets-- one more copy of the street $1$ over
$2$ and one more copy of the street $2$ over $1$ (see Fig 16).
Every positive non-selfintersecting transversal curve will be
expressed in the $m$-dependent basis $a^*,b^*\in \pi_1(T^2_m)$
where $a^*=a^m_1(+), b^*=b^m_1(+)$ as above. Its homology class is
$[\gamma]=ka^*+l(-b^*)\in H_1(T^*_m,Z)$, k>0,l>0. For our
foliation the cycles $a=a^*,b^*=b^{-1}$ are the basic positive
closed transversal curves.

 Every non-selfintersecting positive closed transversal
curve $\gamma\subset T^*_m=T^2 minus (s)$ crosses this fundamental
domain several times (exactly $k+l$ times) from the right to the
left. We denote these segments by $t_j,j=1,...,k+l$ counting them
from the segment $s^+_1=s$ up, along the longest street $0$. Its
canonical representative can be chosen such that it crosses the
domain $2$ exactly $k$ times, entering the domain $0$ in the
points $y_1,...,y_k$; they are naturally ordered by height. It
crosses also the domain $1$ exactly $l$ times, entering the domain
$0$ from $1$  in the points $y_{k+1},...,y_{k+l}$,  We  denote the
equivalent points from the left side of fundamental domain by the
same figures with symbol '. The  equivalent points on the left
side of the street $0$  are the points $y'_{k+1},...,y'_{k+l}$, on
the segment separating street $0$ from the street equivalent to
$1$.. The points $y'_1,y'_2,...,y'_k$ are located on the common
boundary segment separating the second (left) copy of the street
$2$ from the street $0$. Therefore the sequence of segments
ordered by height,is following:

The group (I):

$$t_1=[1,(k+1)'],...,t_l=[l,(k+l)']$$

The group (II):

 $$t_{l+1}=[l+1,1'],...,t_{k}=[k,(k-l)']$$

The group (III):

$$t_{k+1}=[k+1,(k-l+1)'],...,t_{k+l}=[k+l,k']$$

We assign to every segment following homotopy class depending on
the group I,II,III:

 $$\phi:t_j\rightarrow b', t_j\in (I)$$

$$\phi:t_j\rightarrow a',t_j\in (II),(III)$$

We define

$$\phi(\gamma)=\phi(t_{q_1})...\phi(t_{q_{k+l}})\in\pi_1(T^*_m)$$
where $\gamma\sim t_{q_1}...t_{q_{k+l}}$ in the natural order
along the curve. As a result of the previous lemmas, we obtain
following

\begin{thm}
The invariant $\phi (\gamma)$ describes  homotopy classes of the
closed non-selfintersecting positive transversal curves in $T^*_m$
as a positive words written in the free group with two
$m$-dependent generators $a'=a^*,b'=(b^*)^{-1}$ whose transversal
measures are smaller than $m$ . In particular, every indivisible
homology class $k[a']+l[b'],k>0,l>0,$ has unique $m$-dependent
non-selfintersecting  positive transversal representative. It
defines a positive word unique up to the cyclic permutation (i.e.
it defines exactly $k+l$ positive words). This word is calculated
above through the sequence of segments $t_j$. Every such word can
be taken as a part of new Transversal Canonical Basis in $T^*_m$.
\end{thm}

As a corollary, we can see that every integer-valued $2\times
2$-matrix $T,\det T=1,$ with positive entries $k,l,p,q\geq 0$
$$a_1=T(a')=ka'+lb',b_1=T(b')=pa'+qb'$$
determines finite number of different   transversal canonical
bases $a_1,b_1\in \pi_1(T^*_m)$. They are  represented by the
curves $a_1,b_1$ crossing each other transversally in one point
and representing the homology classes $T(a'),T(b')$. We should
locate the  segments of both these curves in the domain number $0$
as above. Let the  curves $a_1,b_1$ are represented by the ordered
sequences of pieces $t_1,...t_{k+l}$ and $t'_1,...t'_{p+q}$
correspondingly going from the right to the left side. We require
existence of  one intersection point for the selected pair
$t_i\bigcap t'_j\neq \emptyset$. All other intersections should be
empty. Every such configuration determines transversal canonical
basis $a_1,b_1$. Its equivalence class  is completely determined
by the relative order of segments $t,t'$ taking into account that
$t_i$ crosses $t'_j$ only once, with intersection index equal to
$1$ from the right to the left side (i.e. $t_i$ is ''higher'' than
$t'_j$ from the left side). Starting from the selected point
$t_i\bigcap t'_j$, we apply the procedure described above in the
theorem. It gives us two positive words $A,B$ in the free group
$F_2$. By definition, the map $\hat{T}:a'\rightarrow
A,b'\rightarrow B$ defines a lift from homology to fundamental
group, for every pair of words $A,B$ constructed in that way.
Consider first a {\bf Reducible}  case where the left ends of the
segments $t_i$ and $t'_j$ crossing each other, are located on the
same part of boundary, where the domain $0$ meets the same domain
$1$ or $2$). If it is the domain $1$, we can see that both words
$A$ and $B$ start with the same letter $b'$. We can deform our
crossing point along the cycle $b'$. After this step we are coming
to the conjugated pair $A',B'$ where $b'$ is sent to the end:
$$A=b'\tilde{A},B=b'\tilde{B}\rightarrow A'=\tilde{A}b',B'=\tilde{B}b'$$.
 The same argument we use if the pair $t_i,t'_j$ ends in the
left boundary of $0$ with domain number $2$ replacing $b'$ by
$a'$. After the series of such steps, we  are coming to the case
where $t_i$ ends in $1$,  $t'_j$ ends in $2$. This process cannot
be infinite because the words $A,B$ are not  powers of the same
word. So this process ends. We call  {\bf Irreducible} the case
where the process ended up.

 Example: The pair of words
$A=b'a'b'a'b', B=b'a'$ requires 5 steps to arrive to the
irreducible case $A'''''=A,B'''''=a'b'$.

In the final irreducible state
 the relation:
$$ABA^{-1}B^{-1}=a'b'a'^{-1}b'^{-1}$$
can be easily seen on the plane with periodic set of segments
removed, looking on the 3-street decomposition of the plane.

 {\bf The Semigroup of unimodular $2\times 2$-matrices with
nonnegative integer-valued entries is free, with two generators
$T_1,T_1$} such that
$$T_1(a')=a'+b',T_1(b')=b',T_2(a')=a',T_2(b')=a'+b'$$
Their lifts to the automorphisms of the free group are following:
$$\hat{T}_1(a')=a'b',\hat{T}(b')=b';\hat{T}_2(a')=a',\hat{T}_2(b')=b'a'$$
Obviously, they both preserve the word $a'b'a'^{-1}b'^{-1}$, so
all semigroup of nonnegative unimodular matrices preserves this
word. It is isomorphic to the semigroup of all positive
automorphisms of free group $F_2$ preserving this word. (We
clarified this question with I.Dynnikov. I don't know where it is
written).

Comparing this result with the previous arguments, we are coming
to the following

{\bf Conclusion}: The Semigroup $G$ of all positive automorphisms
of the free group $F_2$ in the alphabet $a',b'$ preserving the
conjugacy class of the word $\kappa=a'b'a'^{-1}b'^{-1}$, contains
following parts:

1.A free semigroup $G_{\kappa}$ consisting of  transformations
preserving this word exactly; It is isomorphic to the semigroup of
matrices $T$ with $\det T=1$ and nonnegative integer entries;

2. Every element    $T\in G_{\kappa}$ defines  finite number of
positive transformations  $T'\in G$ such that the corresponding
pair of positive words $A',B'$ of $T'\in G$ are simultaneously
conjugate to the words $A,B$ of $T\in G_{\kappa}$. All these
conjugations should be performed by the simultaneous cyclic
permutation of the words $A,B$ removing the same letter from the
end  and sending it to the beginning (until both words ends with
the same letter). Total number of these elements (conjugations) is
equal to the sum of matrix elements of $T$ minus 2. The natural
projection of semigroup $G$ into $G_{\kappa}$ is such that the
inverse image of each nonnegative unimodular matrix  $T$ contains
exactly $k+l+p+q-2$ positive automorphisms. In the example above
we have $k+l+p+q-2=5$.

 \vspace{0.3cm}

{\Large Appendix. Existence of the Transversal Canonical Basis}

\vspace{0.2cm}

Consider any  $C^1$-smooth vector field on the compact smooth
Riemann Surface $M_g$ of the genus $g\geq 1$, with saddle
(nondegenerate) critical points only. We introduce  {\bf the Class
$G$ of vector fields} requiring that there is no saddle
connections and no periodic trajectories for the vector fields in
this class. Let a finite family of smooth disjoint
non-selfintersecting positively oriented segments $s_j\subset M_g$
is given, $j=1,...,N$, transversal to the vector field in every
point, including the boundary points $P_{1,j}\bigcup
P_{2,j}=\partial s_j$. We think that the points $P_{1,j}$ are left
boundaries of the segments $s_j$.

\begin{lem}
There exists a non-selfintersecting smooth closed curve $\gamma$
transversal to our vector field, which does not  cross any segment
$s_j$.
\end{lem}
Proof. Start the trajectory $\kappa$  anywhere outside of the
segments $\bigcup s_j$. It either meets first time some segment
$s_j$ of that family at the moment $t<\infty$, or never meets this
family for $t\rightarrow\infty$. In the second case we construct a
closed transversal curve $\gamma$ as usually, ignoring the
segments $s_{j_1}$. In the first case we make a turn
contr-clockwise and move from the point $\kappa(t-\epsilon)$
towards the end $P_{1,j_1}$, parallel to the segment $s_{j_1}$,
constructing the almost transversal curve $\kappa$. After passing
a small transversal distance $\delta$ beyond the point $P_{1,j_1}$
to the left, we make a second turn and begin to move straight
ahead extending the almost transversal curve $\kappa$ along the
trajectory of vector field. The small constants $\epsilon, \delta$
are chosen a priori. They remain unchanged during this
construction: they are the same for all steps. At some moment we
either reach the segment $s_{j_2}$ or never meet this family. The
second case was already taken into account. In the first case we
repeat our move with the same parameters $\delta,\epsilon$: we
make  turn and move beyond the left end. After $\delta$-passing
it, we make one more turn and extend our transversal curve along
the trajectory, and so on. The infinite almost transversal curve
$\kappa$ constructed by this process, certainly has an
$\omega^+$-limiting point $P\in M_g$. It cannot be an inner point
of any segment $s_k$. There are 3 possibilities:

The case one: $P$ is a left end of some segment $s_k$.

The case two: $P$ is a right end of some segment $s_k$.

The case three: $P$ is some point distant from $\bigcup s_k$.

The first case is trivial: There is a sequence of points
$\kappa(t_q)\rightarrow P$ for $t_q\rightarrow \infty$ which
approach  $P$ from the left side of it. Therefore the small
transversal intervals $r_q=[\kappa(t_{q+1})\kappa(t_q)]$ are
positive. The  curve $\kappa(t),t_q\leq t\leq t_{q+1},r_q$ is
almost transversal, positive and closed.

In the second case our curve $\kappa$ approaches $P$ from the
right side by sequence of points located inside of the pieces of
trajectories of vector field. Assuming that these pieces are not
infinite for $t\rightarrow -\infty$, we conclude that they already
met some segment of our family at negative moment of time.
Therefore we see that all of them are coming after the
$\delta$-passing some of the left ends of our segments. Therefore
they start to cross each other already at some  finite time
$t'<\infty$. The curve $\kappa$ is selfintersecting. Taking this
curve till the first self-crossing, we obtain a closed almost
transversal curve $\gamma$.

The third case is either trivial if trajectory passing the point
$P$ never meets our segments, or reduces to the cases one and two.

 Lemma is proved.

Now we cut the surface $M_g$ along the transversal closed curve
$\gamma\rightarrow a^+\bigcup a^-$.  We obtain a surface
$\bar{M_g}$ such that $\partial \bar{M_g}=a^+\bigcup a^-$ as in
the Fig 1. The trajectories start at the in-cycle $a^+$ and end up
in the out-cycle $a^-$, but there is an obstacle $s$ inside.

\begin{lem}
There exists a closed curve $\gamma_1\subset M_g$ transversal to
foliation such that it crosses $\gamma$ transversally in one point
only; it does not cross the segment $s$.

\end{lem}

Proof. We are going to construct a trajectory which starts and
ends up in the boundary $\partial\bar{M_g}=a^+\bigcup a^-$ and
does not touch the segment $s$. This construction  immediately
implies the existence of the transversal curve $\gamma_1$. Let
every trajectory  started in $a^+$, always meets $s$. Find the
first point which trajectory meets and move it to the left  along
the transversal segment $s$ as far as possible. Either the
trajectory joining us with $a^+$ meets the second end $P_2$, or we
can continue moving left. Every time we jump (if needed) to the
nearest point to the boundary $a^+$ if new point appears. Finally
we are coming to conclusion that either the most left point
$P_1\in s$ can be reached directly from the boundary $a^+$, or the
most right $P_2$. We denote this point $P_k$. Changing time
direction, we prove that we can directly reach the second boundary
$a^-$ from another boundary point $P_l\in\partial s,l\neq k$, not
crossing $s$. So we construct a transversal connection from $a^+$
to $a^-$ moving from $a^+$ to $P_k$ along the trajectory found
above, after that we move along the segment $s$ (not touching it),
and finally we move from the point near $P_l$ to $a^-$ along the
second trajectory found above. Here $k,l=1,2$. We make this curve
closed in $M_{g-1}$ using additional transversal piece along the
cycle $a^-$ of the desired orientation.

 Lemma is proved.

Let us prove following

\begin{thm}
For every vector field of the class $G$ on the Riemann Surface
$M_g$ with genus $g\geq 3$, there exists an incomplete Transversal
Canonical Basis $a_1,b_1,a_2,b_2,a$ such that all basic cycles are
non-selfintersecting and transversal to the vector field. All
pairwise intersections are transversal and nonempty only for
$a_k\bigcap b_k$ where they consist of one point, $k=1,2$.
\end{thm}
Proof.

 Step 1: Reduction of genus. We construct a first pair
of the transversal basic cycles $a_1,b_1$. After that we cut $M_g$
along these cycles $a_1\bigcup b_1$. This operation leads to the
manifold $W$ with boundary $A_1=A_1^+\bigcup A_1^-=\partial W$
consisting of parallelogram in the plane $C_1$ as before.  We
construct the transversal segments $s^{\pm}$ joining the pair of
saddles as in the Fig 5. Now we remove all part of $W$ outside of
the segments $s^{\pm}$. We are coming to the manifold $$\bar{W}
,\partial\bar{W}=s^+\bigcup s^-$$ Finally, we identify the
segments $s^+$ and $s^-$ and obtain from $\bar{W}$ a closed
2-manifold $M_{g-1}$ with vector field inherited from $M_g$ (which
easily can be done $C^1$), and transversal segment $s^+=s^-\subset
M_{g-1}$

Step 2: Construction of the second transversal pair $a_2,b_2$. We
use for that the lemmas above. They explain how to deal with the
case of one transversal segment $s$. After that we reduce  genus
second time and come to the manifold $M_{g-2}$ with $C^1$-smooth
vector field and  two transversal segments $s_1,s_2\subset
M_{g-2}$

Step 3: Construction of the transversal cycle $a\subset M_{g-2}$
not crossing the previous segments. We use for that the lemma
above.

Our Theorem follows from these steps.

Combining this theorem with lemmas 5 and 6 above (see  Section 2),
we obtain the following result for $g=3$.

\begin{thm}
For every Hamiltonian System on the Riemann Surface given by the
real part of the generic holomorphic one-form there exists a
Transversal Canonical Basis.
\end{thm}

{\bf As G.Levitt pointed out to the present author, his work
\cite{L} published in 1982, can be used to construct Transversal
Canonical Bases on every Riemann Surfaces of any genus $g$ for the
Dynamical Systems of the class $G$ (above)}. Let us present here
this construction.
\begin{thm} (G.Levitt, private communication) For every generic dynamical system
of the class $G$ there exist a Transversal Canonical Basis
\end{thm}

 Proof. The result of \cite{L} is following: For every dynamical system of this class
there exist exactly $3g-3$ non-selfintersecting and pairwise
non-intersecting  transversal cycles $A_i,B_i,C_i, i=1,...,g-1,$
where $C_{g-1}=C_1$. These cycle bound two sets of surfaces:
$A_i\bigcup B_i\bigcup C_i=\partial P_i$ where $P_i$ is a genus
zero surface (''pants''). The trajectories enter pants $P_i$
through the cycles $A_i,B_i$ and leave it through the cycle $C_i$.
There is exactly one saddle inside of $P_i$ for each number $i$.
The cycles $A_i,B_i,C_{i-1}$ bound also another set of pants $Q_i$
where $\partial Q_i=A_i,B_i,C_{i-1}$. The trajectories enter $Q_i$
through $C_{i-1}$ and leave it through $A_i$ and $B_i$. There is
also exactly one saddle inside $Q_i$.

The construction of Transversal Canonical Basis based on that
result is following: Define  $a_g$ as $a_g=C_1$. Choose a segment
of trajectory $\gamma$ starting and ending in $C_1$ assuming that
it passed all this ''necklace'' through $P_k,Q_k$. For each $i$
this segment meets either $A_i$ or $B_i$. Let it meets $B_i$ (it
does not matter). We define $a_i$ as $a_i=A_i,i=1,...,g-1$. Now we
define $b_g$ as a piece of trajectory $\gamma$ properly closed
around the cycle $C_1=a_g$. This curve is almost transversal in
our terminology (above), so its natural small approximation is
transversal. We are going to construct the cycles $b_i$ for
$i=1,...,g-1$ in the union $P_i\bigcup Q_i$. Consider the saddle
$q_i\in Q_i$. There are two separatrices leaving $q_i$ and coming
to $A_i$ and $B_i$ correspondingly. Continue them until they reach
$C_i$ through $P_i$. Denote these pieces of separatrices by
$\gamma_{1,i},\gamma_{2,i}$. Find the segment $S_i$ on the cycle
$C_i$ not crossing the curve $\gamma$ chosen above for the
construction of the cycle $b_g$. The curve
$\gamma_{1,i}S_i\gamma_{2,i}^{-1}$ is closed and transversal
everywhere except the saddle point $q_i$. We approximate now the
separatrices $\gamma_{1,i}$ and $\gamma_{2,i}^{-1}$ by the two
pieces of nonseparatrix trajectories $\gamma'_i,\gamma''_i$
starting nearby of the saddle $q_i$ from one side of the curve
$\gamma_{1,i}\gamma_{2,i}^{-1}$. Close this pair by the small
transversal piece $s_i$. There are two possibilities here (two
sides). We choose the side such that the orientation of the
transversal piece $s_i$ near the saddle $q_i$ is agreed with the
orientation of the segment $S_i$, so the whole curve $b_i=s_i
\gamma'_iS_i\gamma''_i$ is closed and almost transversal. We
define $b_i$  as a proper small transversal approximation of that
curve. So our theorem is proved because every cycle $b_i$ crosses
exactly one cycle $a_i$ for $i=1,2,...,g$.

\begin{center}
\mbox{\epsfxsize=5cm \epsffile{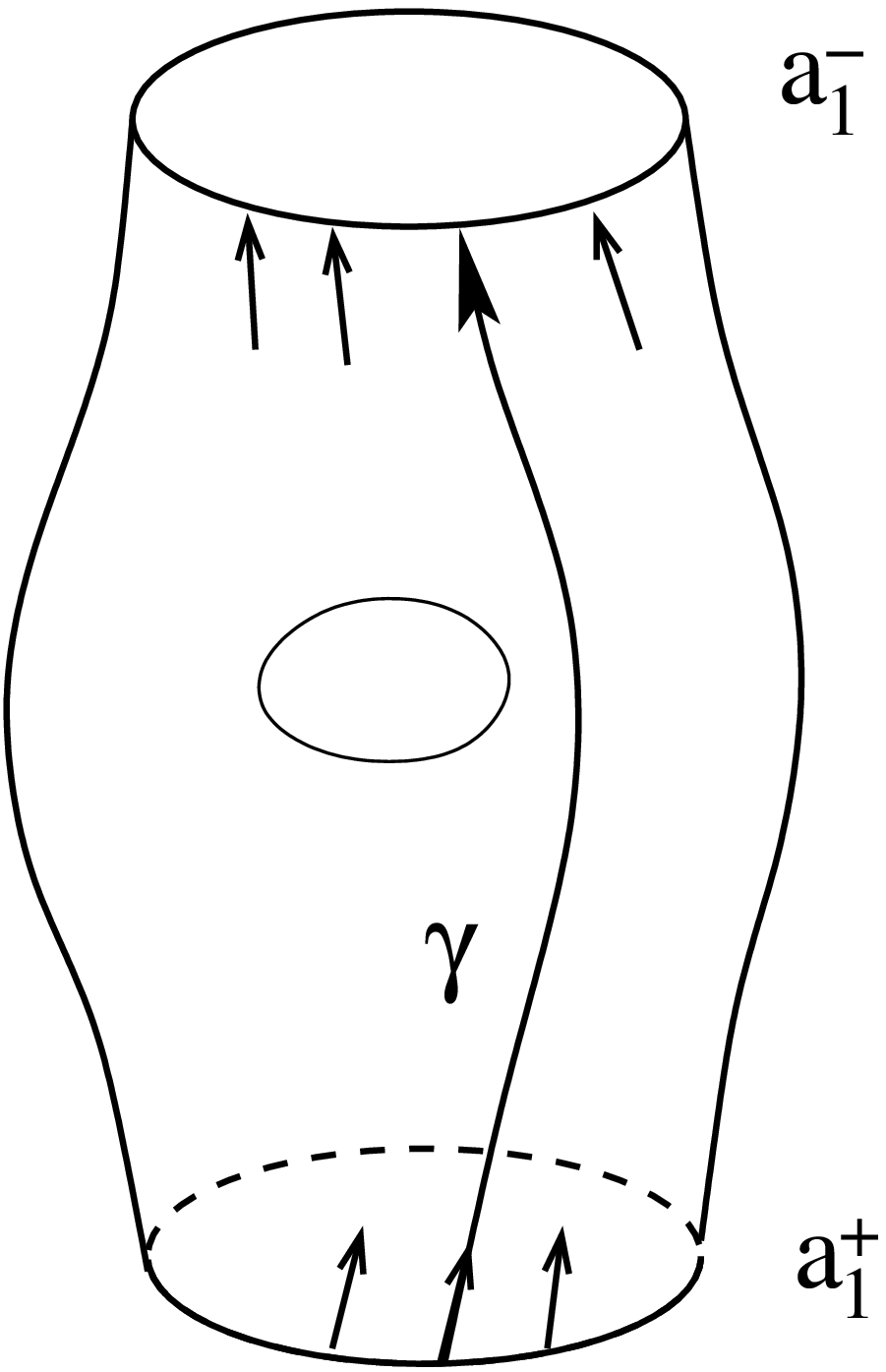}}

Fig 1.
\end{center}

\begin{center}
\mbox{\epsfxsize=4cm \epsffile{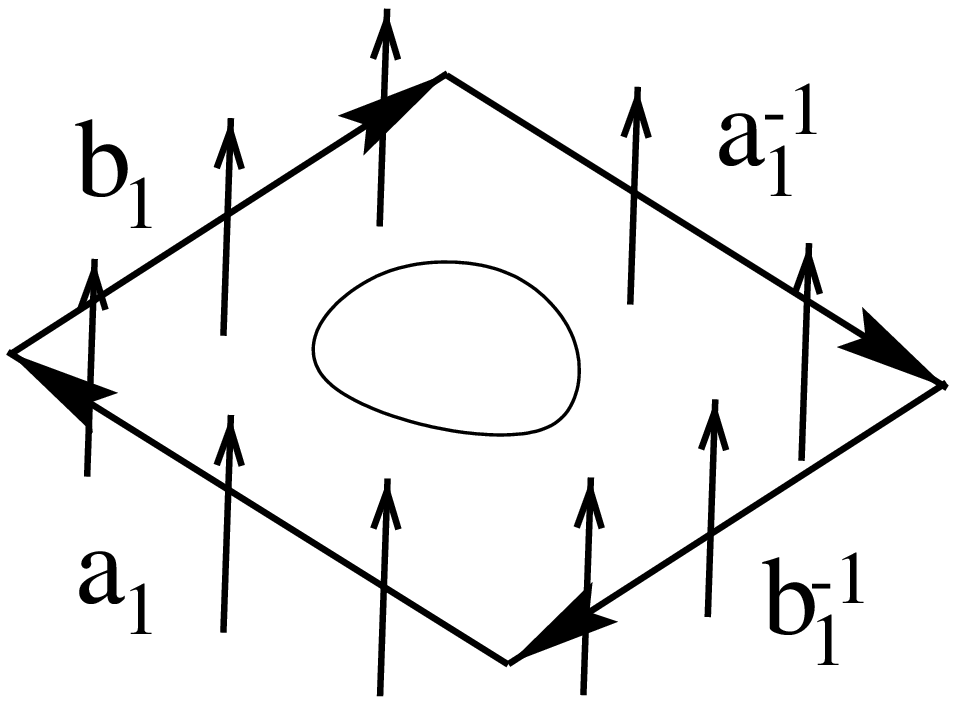}}

Fig 2.
\end{center}

\begin{center}
\mbox{\epsfxsize=10cm \epsffile{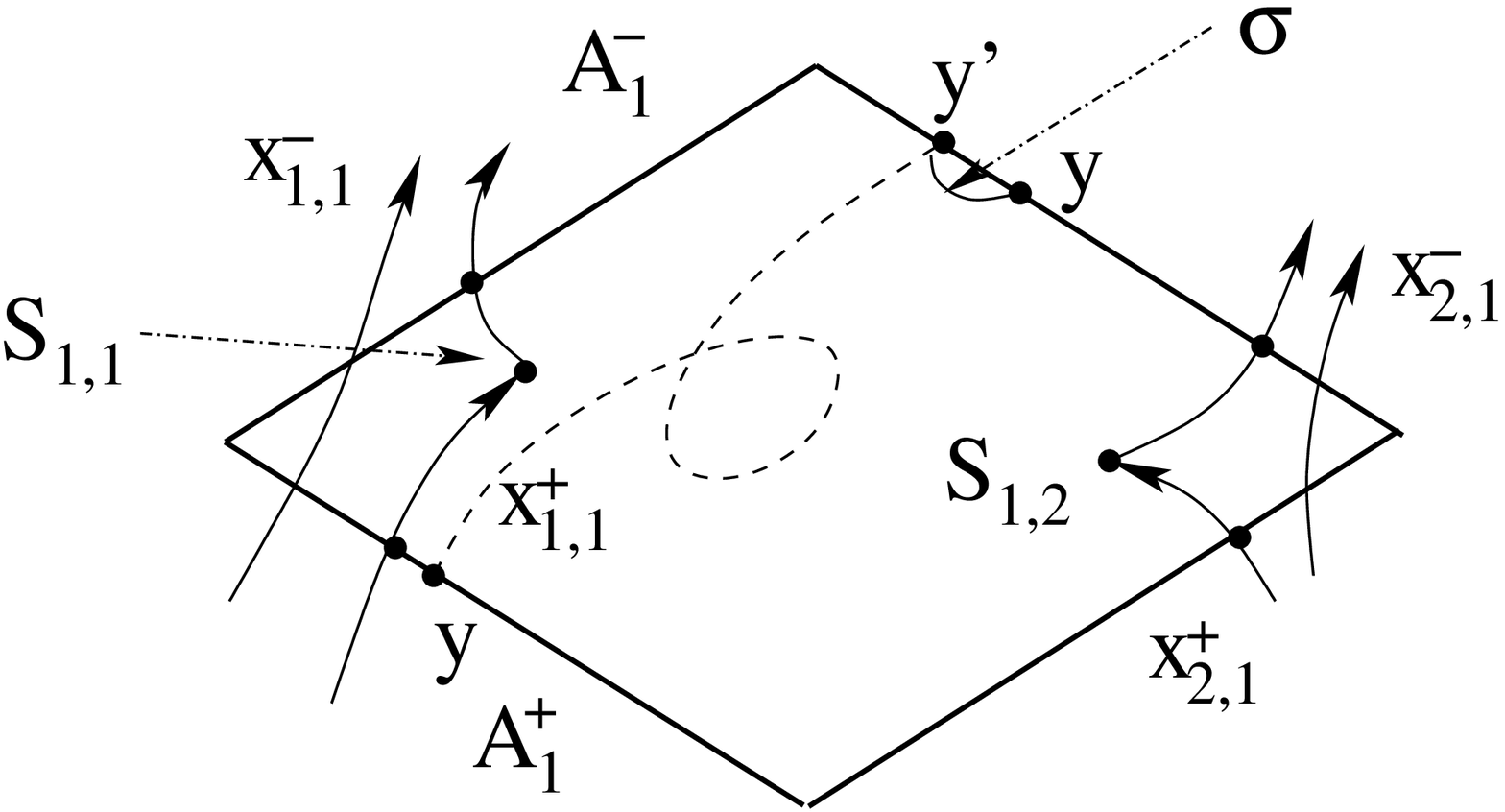}}

Fig 3.
\end{center}

\begin{center}
\mbox{\epsfxsize=10cm \epsffile{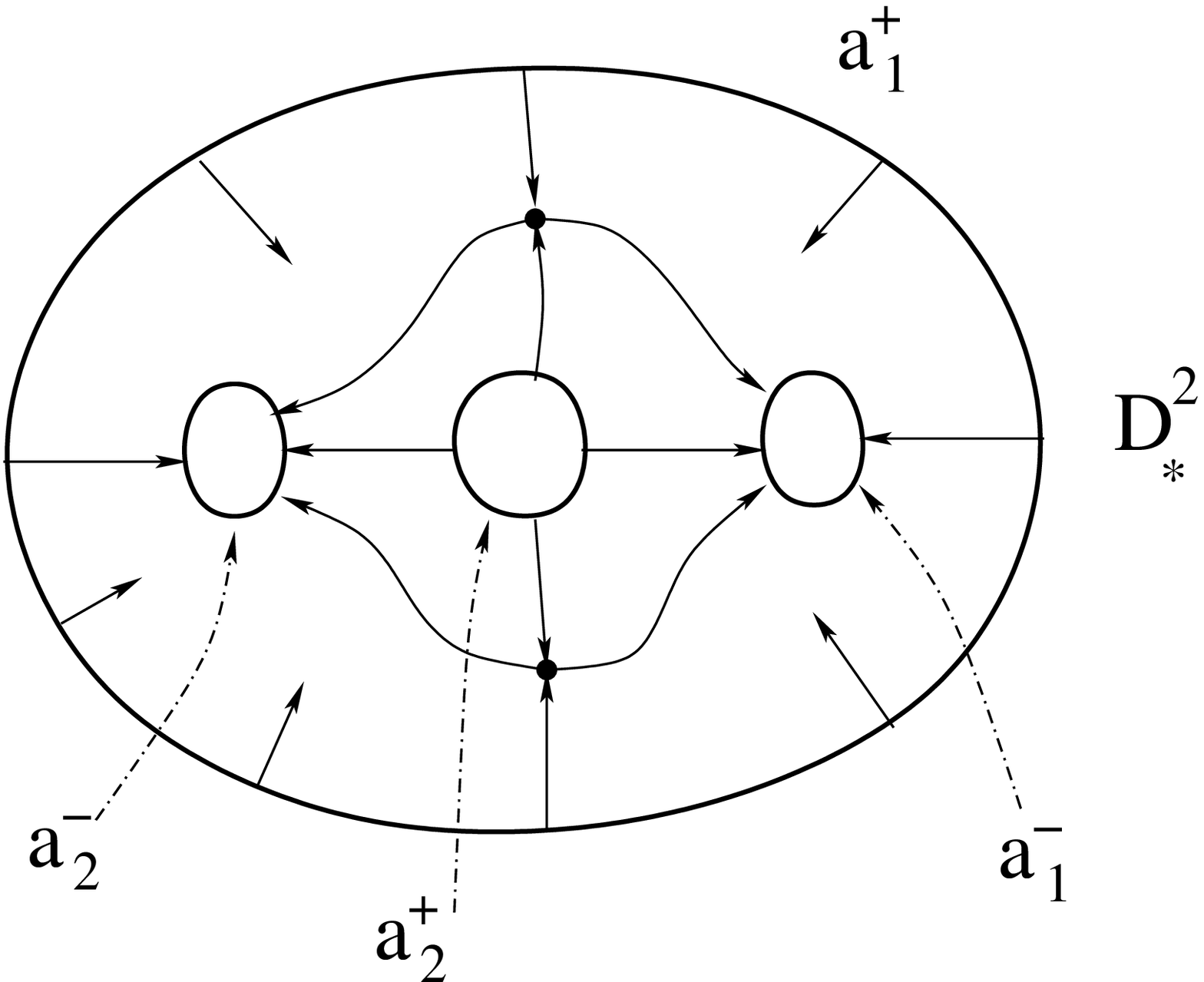}}

Trajectory connections: $a_1^+\rightarrow a_1^-$,
$a_2^+\rightarrow a_2^-$.

Fig 4a.
\end{center}

\begin{center}
\mbox{\epsfxsize=10cm \epsffile{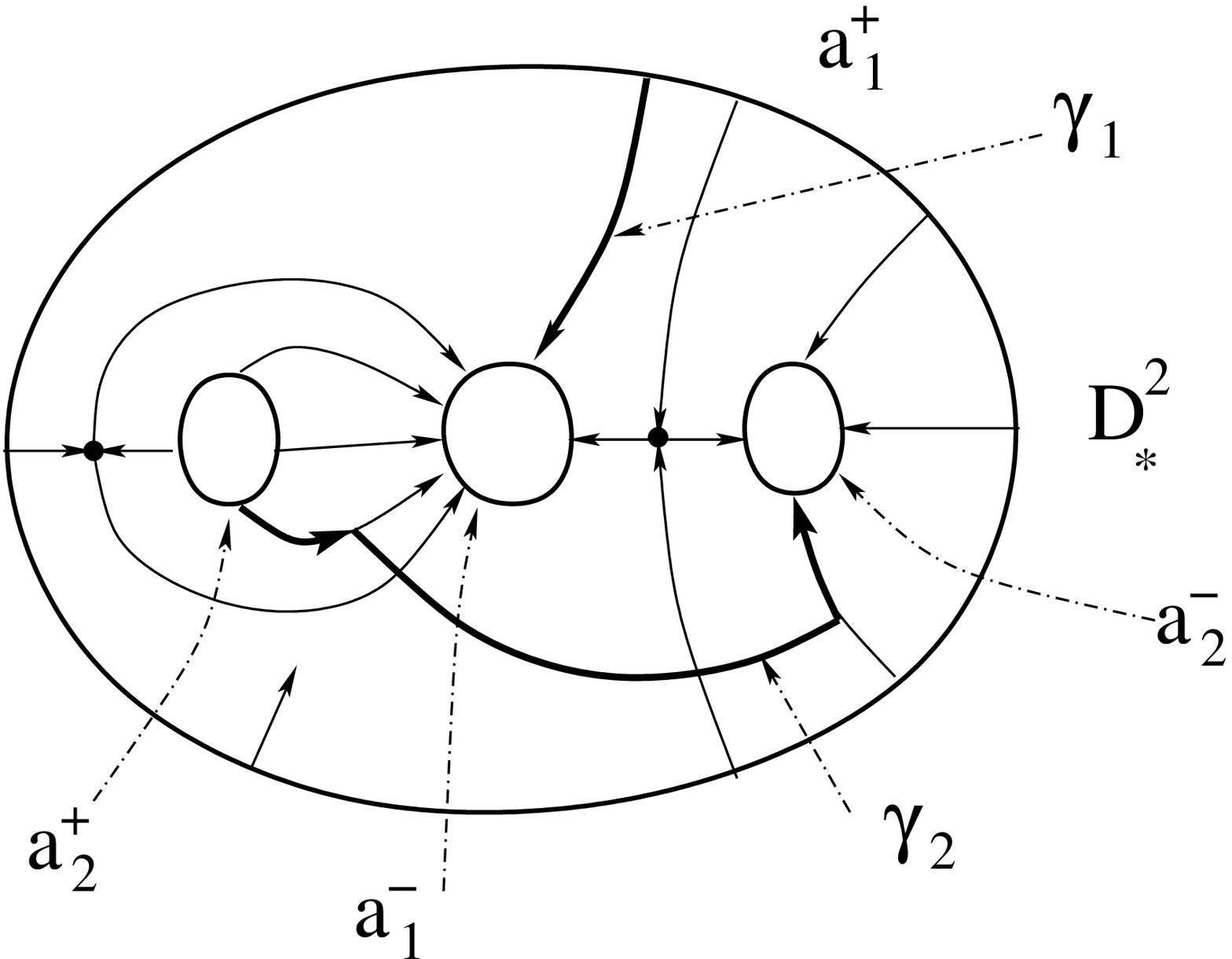}}

Trajectory connection: $a_1^+\rightarrow a_1^-$.

Almost transversal curve $\gamma$ (boldface): $a_2^+\rightarrow
a_2^-$.

Fig 4b.
\end{center}

\begin{center}
\mbox{\epsfxsize=8cm \epsffile{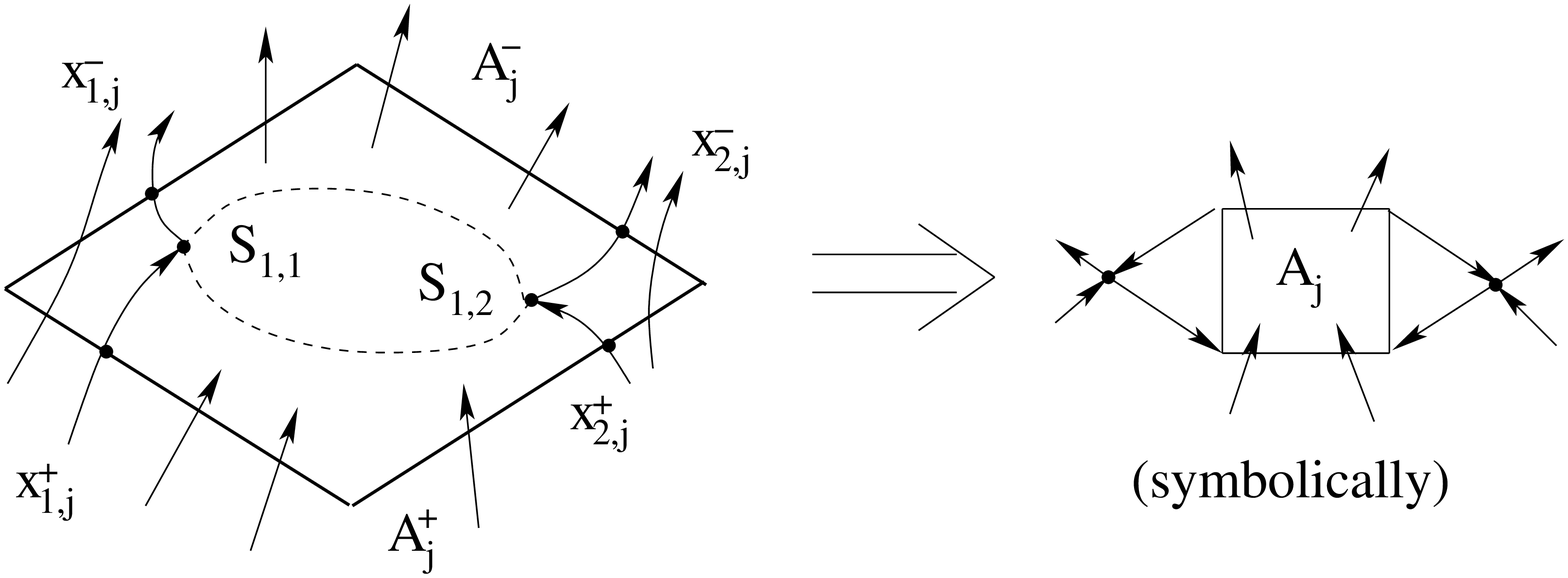}}

Fig 5.
\end{center}

\begin{center}
\mbox{\epsfxsize=7cm \epsffile{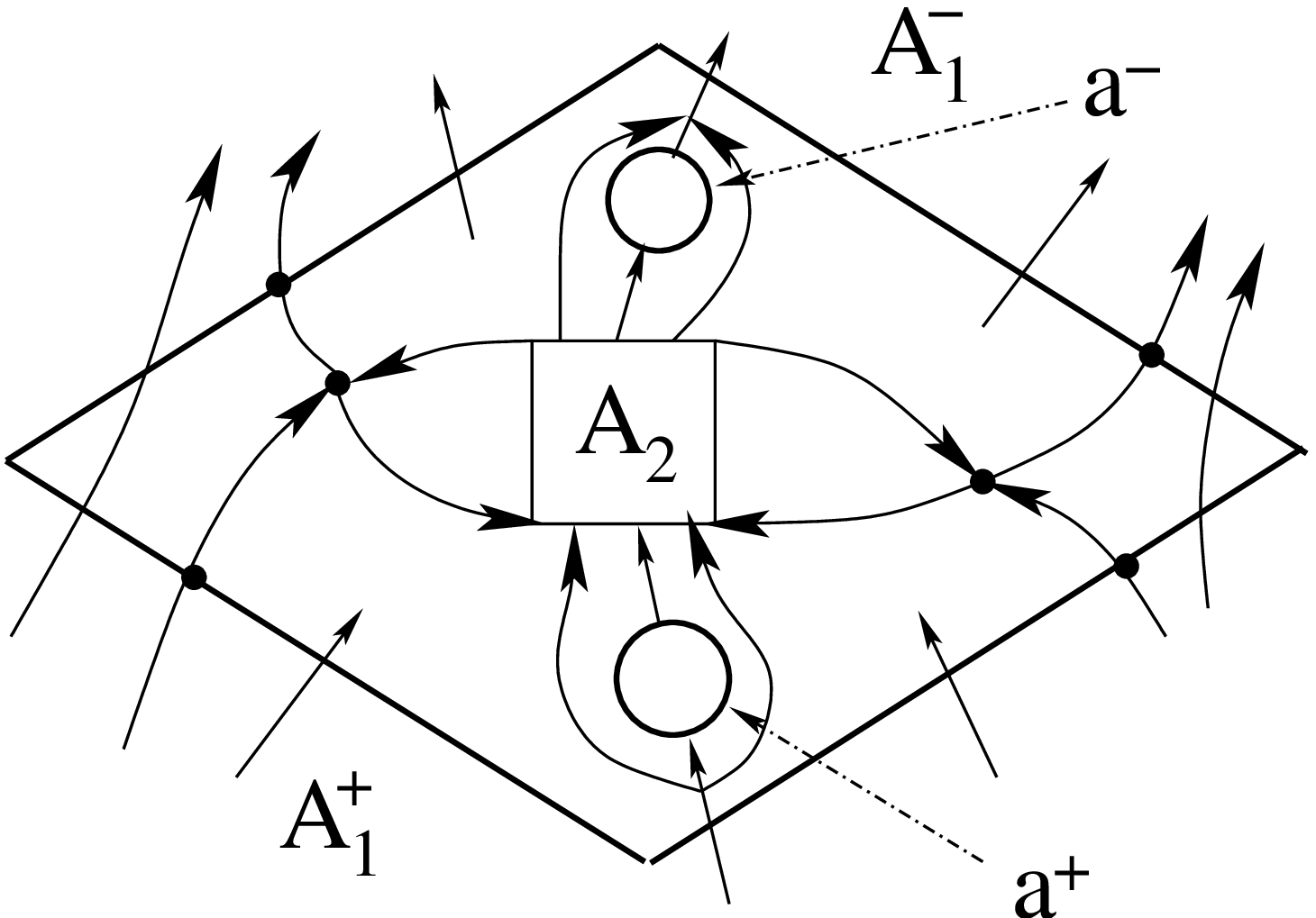}}

Non-extendable diagram of the type $T^2$.

Fig 6.
\end{center}

\begin{center}
\mbox{\epsfxsize=9cm \epsffile{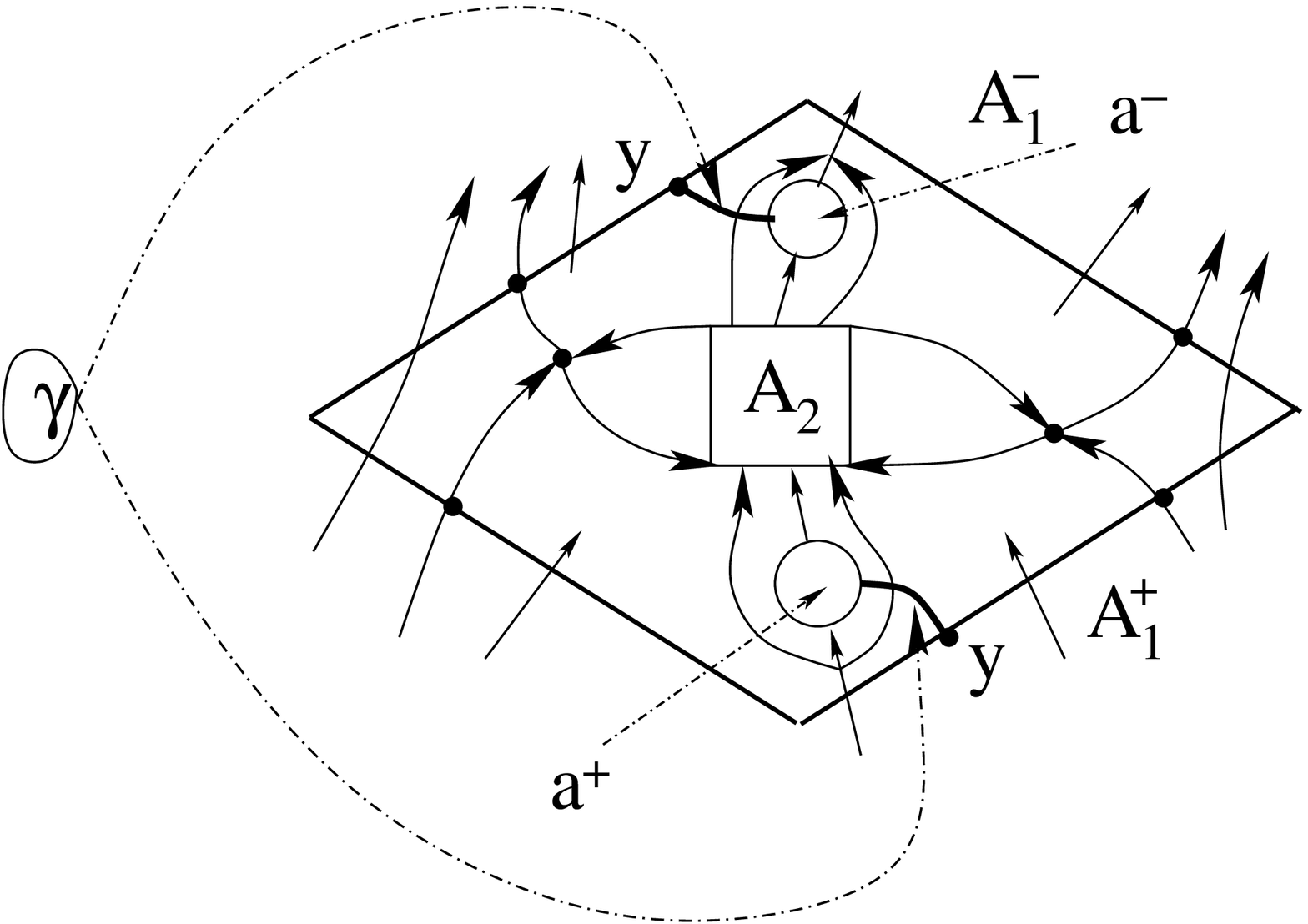}}

Reconstruction: new transversal curve $\gamma:$ $a_1^+\rightarrow
a_1^-$ is passing through the boundary.

Fig 7.
\end{center}

\pagebreak

\begin{center}
\mbox{\epsfxsize=6cm \epsffile{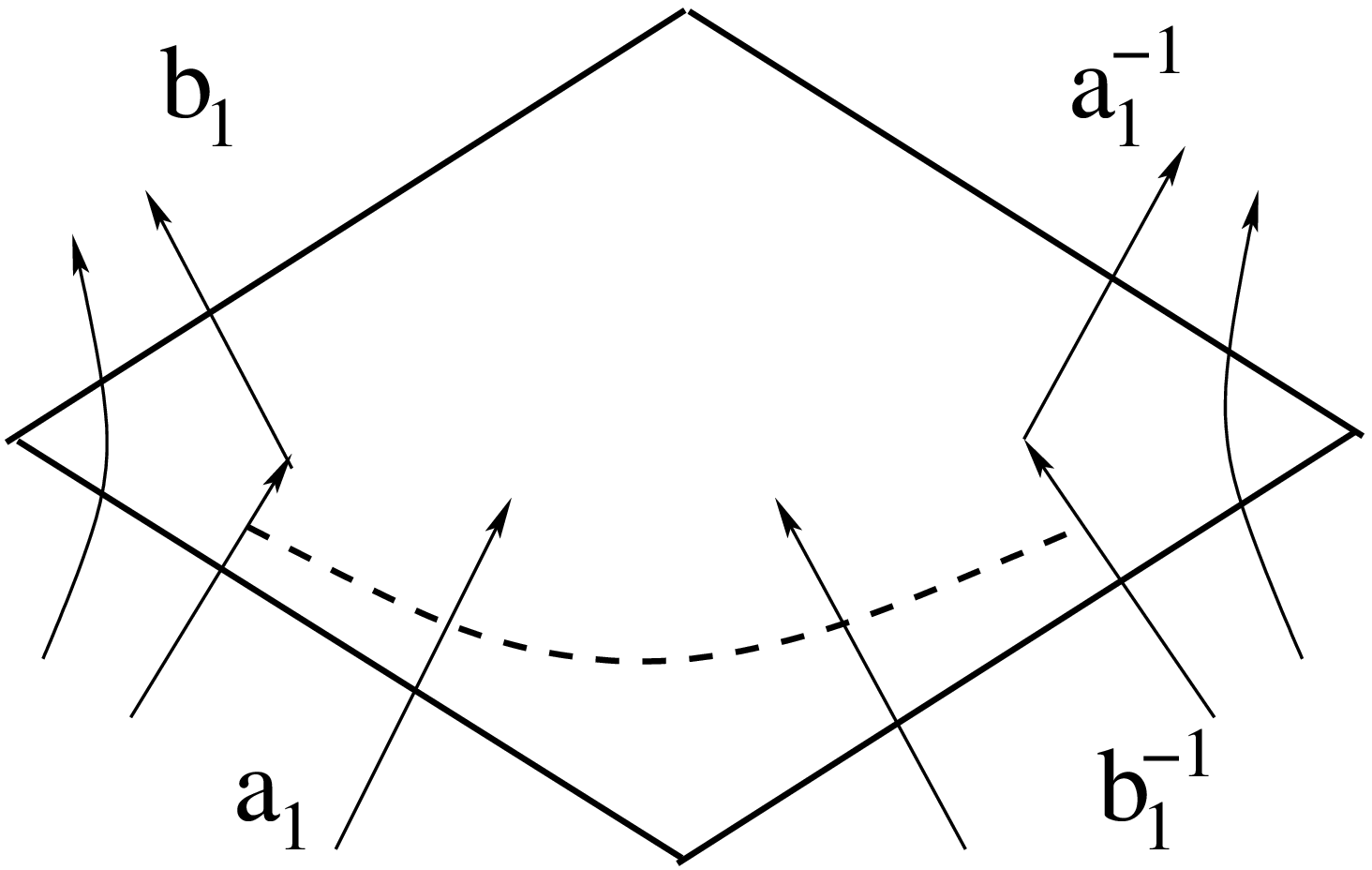}}

Fig 8.
\end{center}

$g=2$ 

\vspace{-1cm}
\begin{center}
\mbox{\epsfxsize=14cm \epsffile{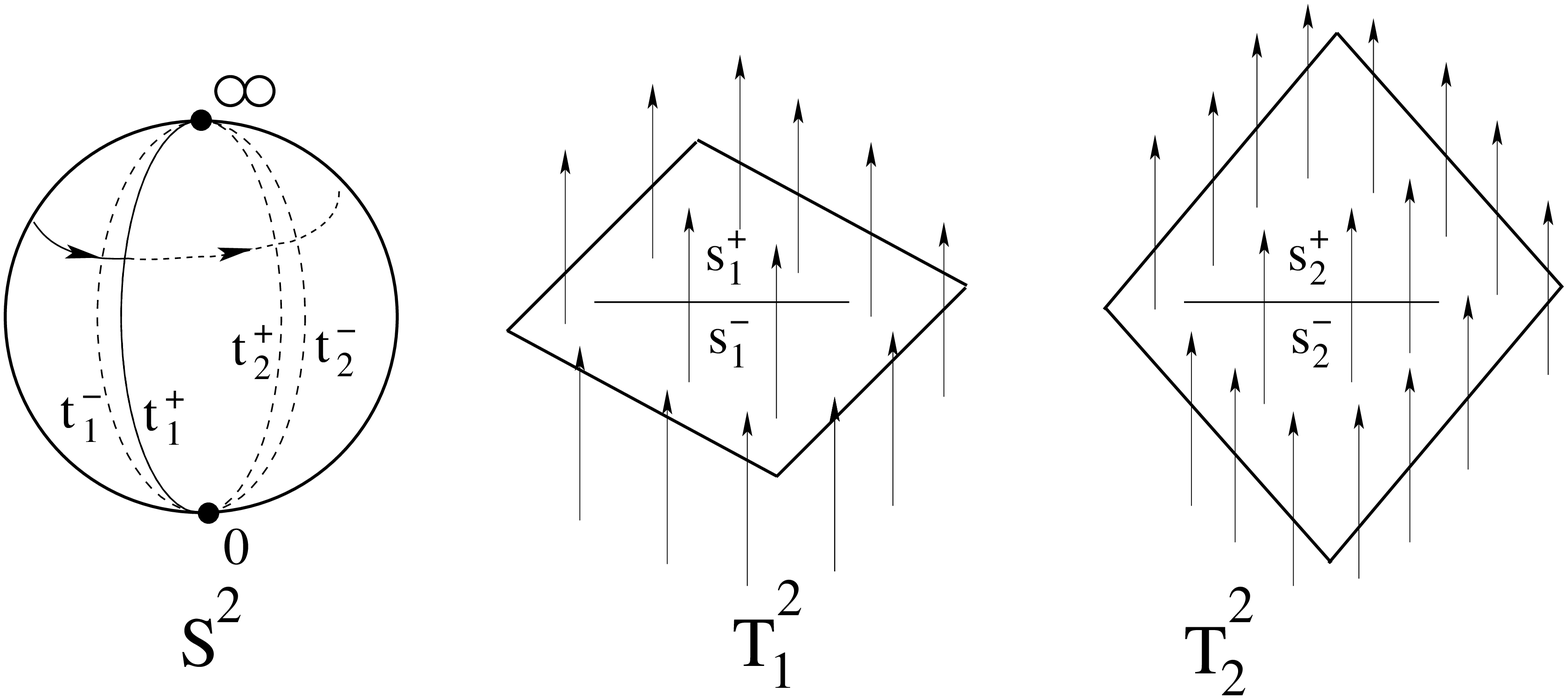}}

(Sperical parts can be dropped)

Fig 9a.
\end{center}

$g=3$

\vspace{-1cm}
\begin{center}
\mbox{\epsfxsize=9cm \epsffile{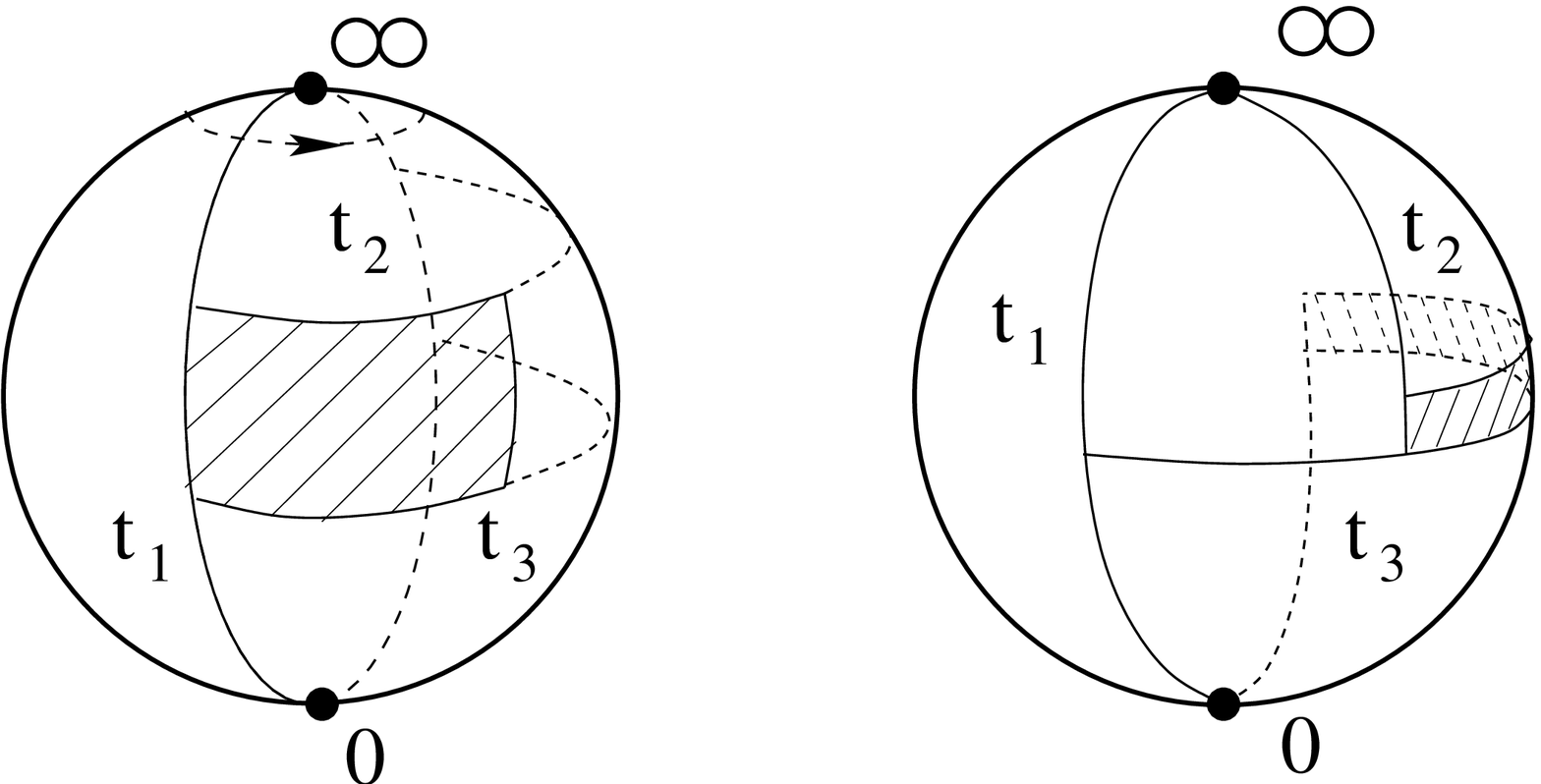}}

``corridors'' $t^+_i\rightarrow t^-_i$.

Fig 9b.
\end{center}

\pagebreak

Graphs

\vspace{-0.5cm}
\begin{center}
\mbox{\epsfxsize=9cm \epsffile{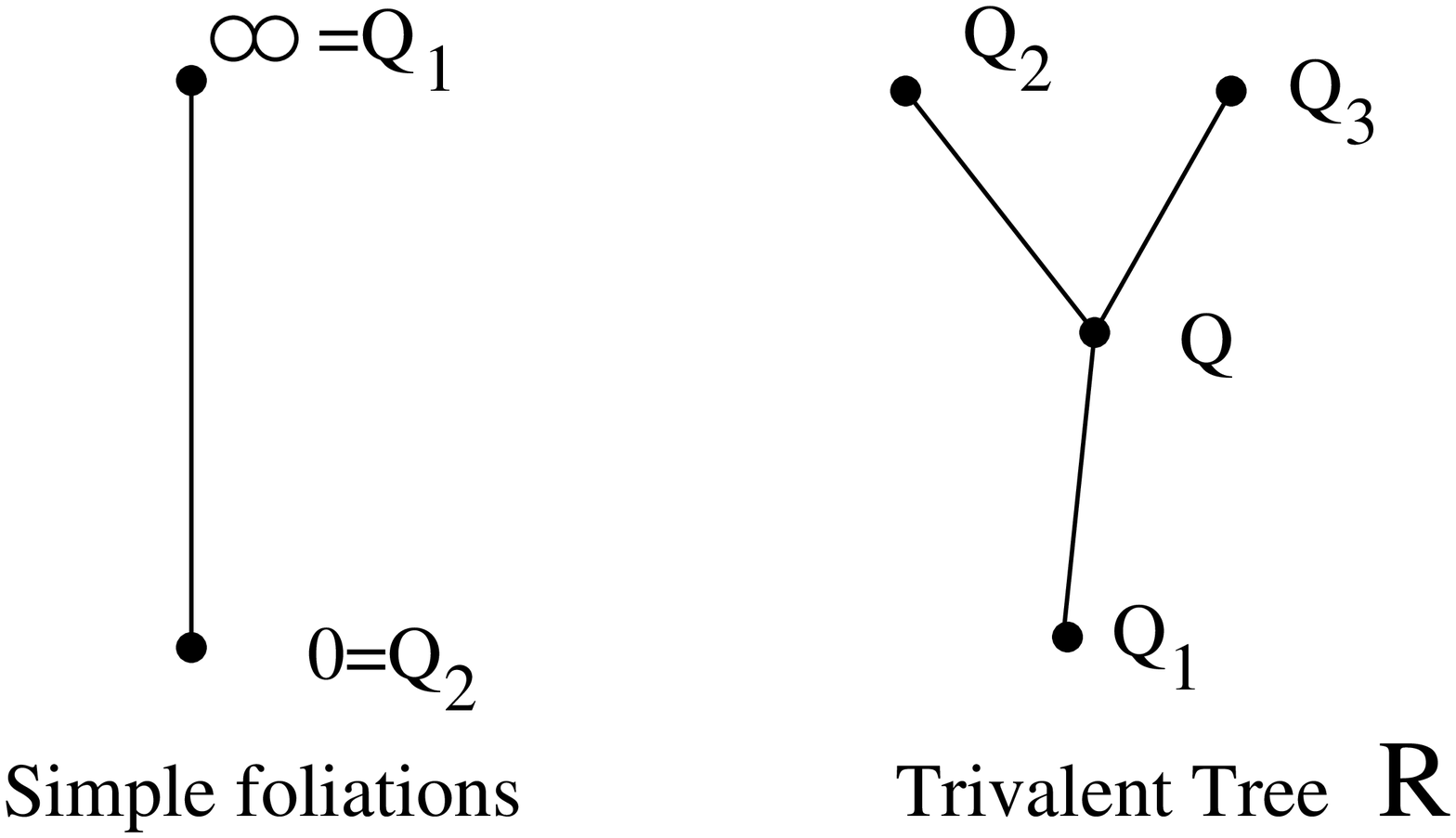}}

Fig 9c.
\end{center}

\begin{center}
\mbox{\epsfxsize=6cm \epsffile{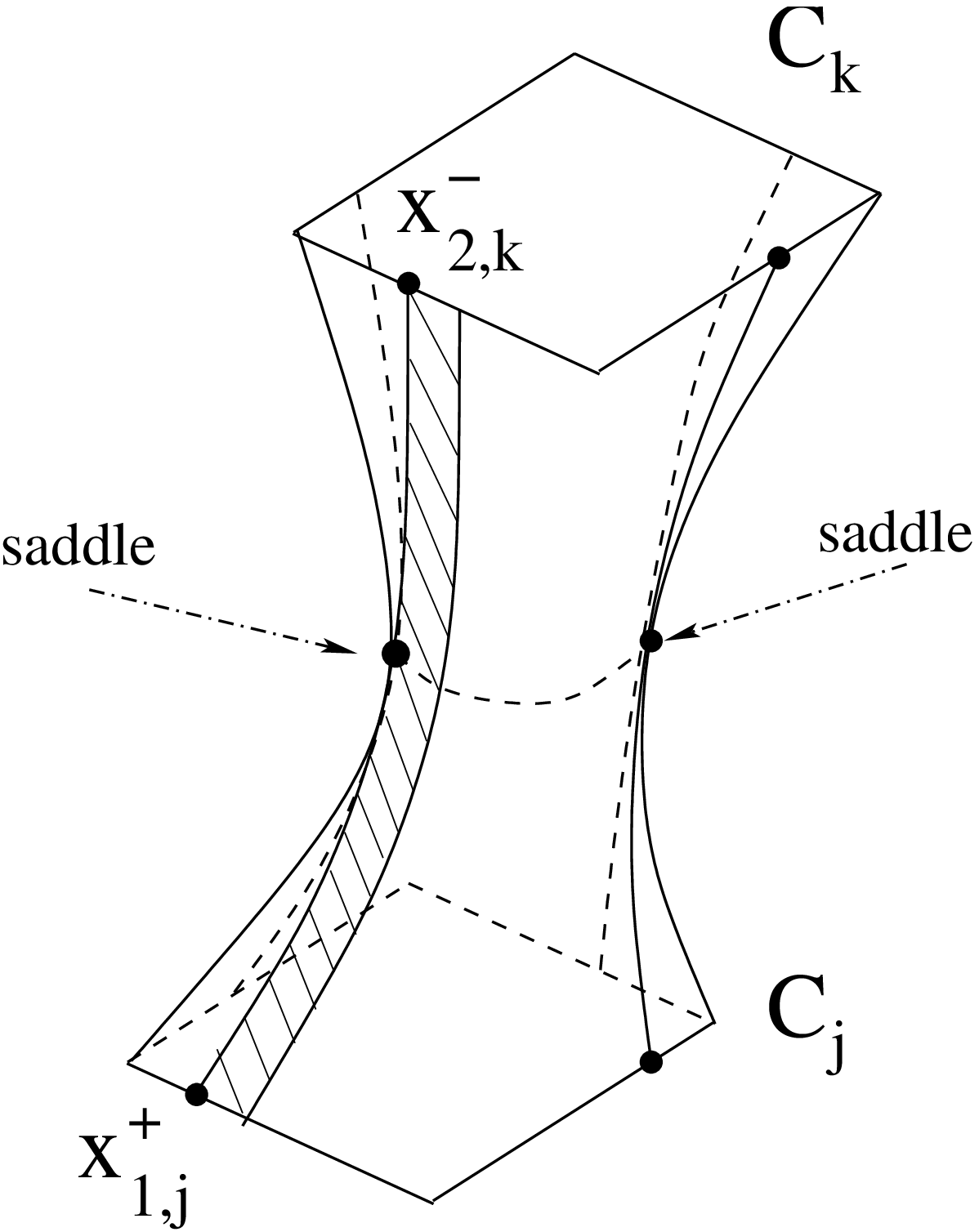}}

Fig 10.
\end{center}

\begin{center}
\mbox{\epsfxsize=12cm \epsffile{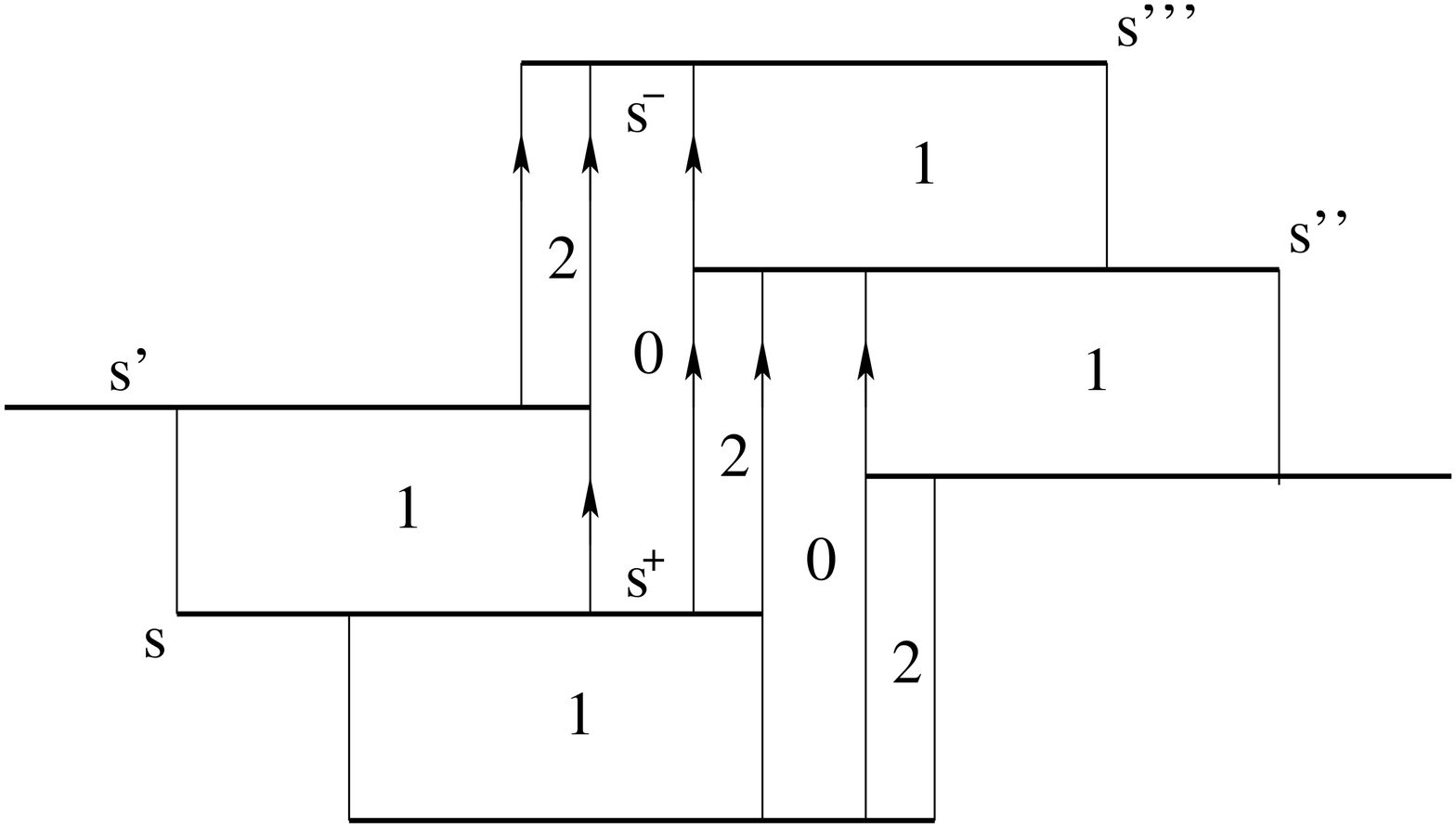}}

The 3-street model.

Fig 11.
\end{center}

\begin{center}
\mbox{\epsfxsize=5cm \epsffile{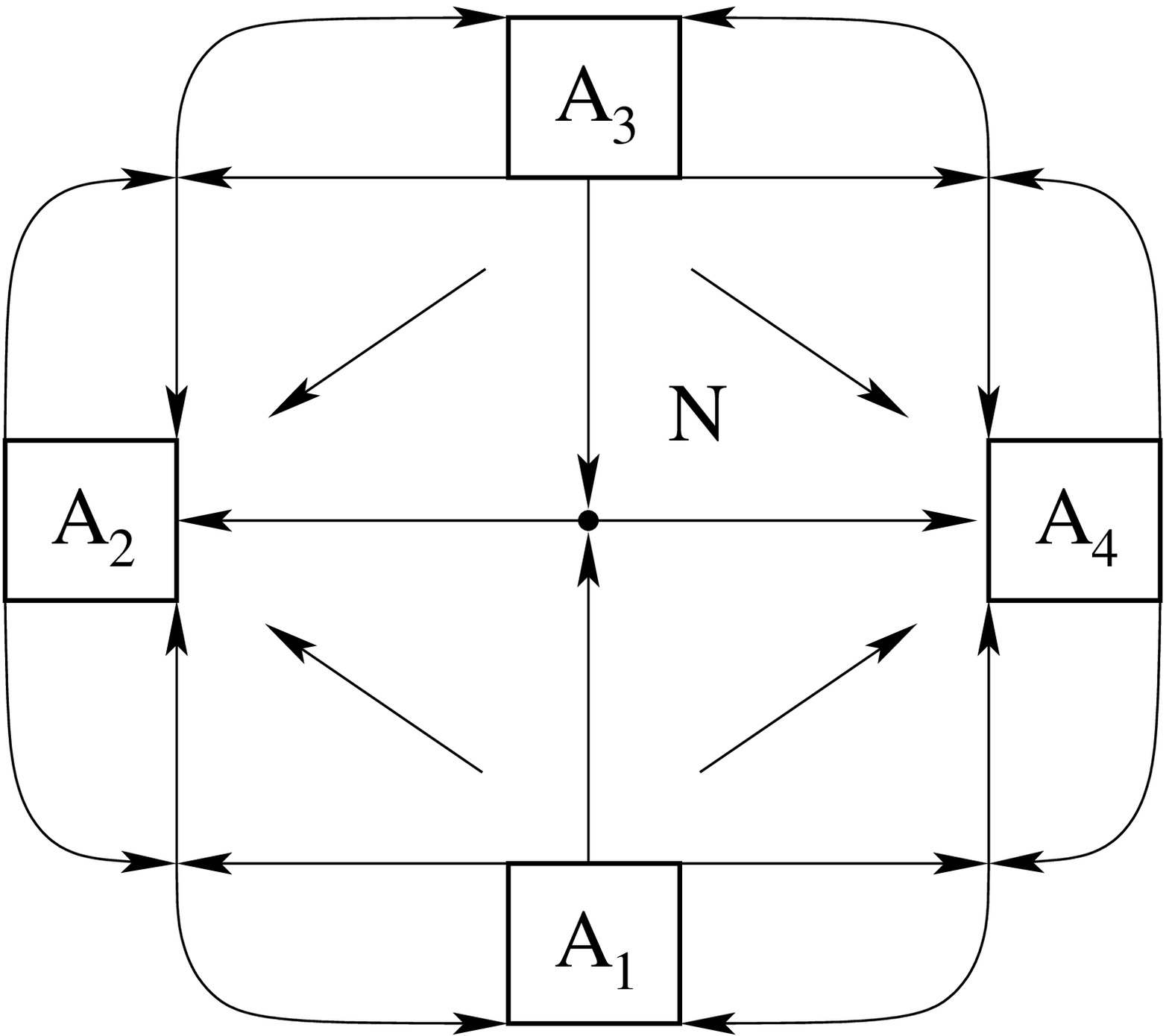}}\hspace{1cm}
\mbox{\epsfxsize=5cm \epsffile{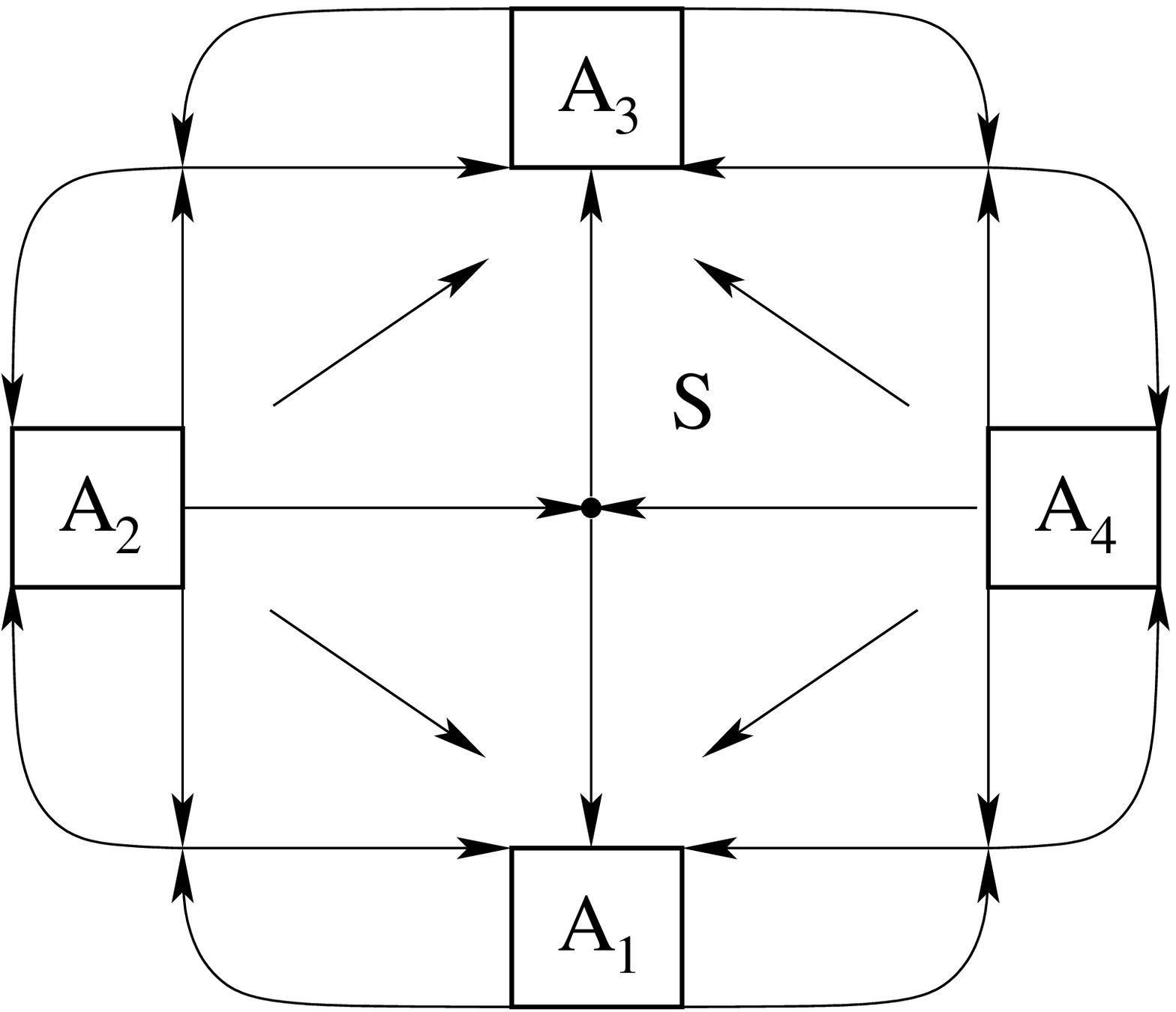}}

North hemisphere \hspace{4cm} South hemisphere

Fig 12.
\end{center}

\begin{center}
\mbox{\epsfxsize=9cm \epsffile{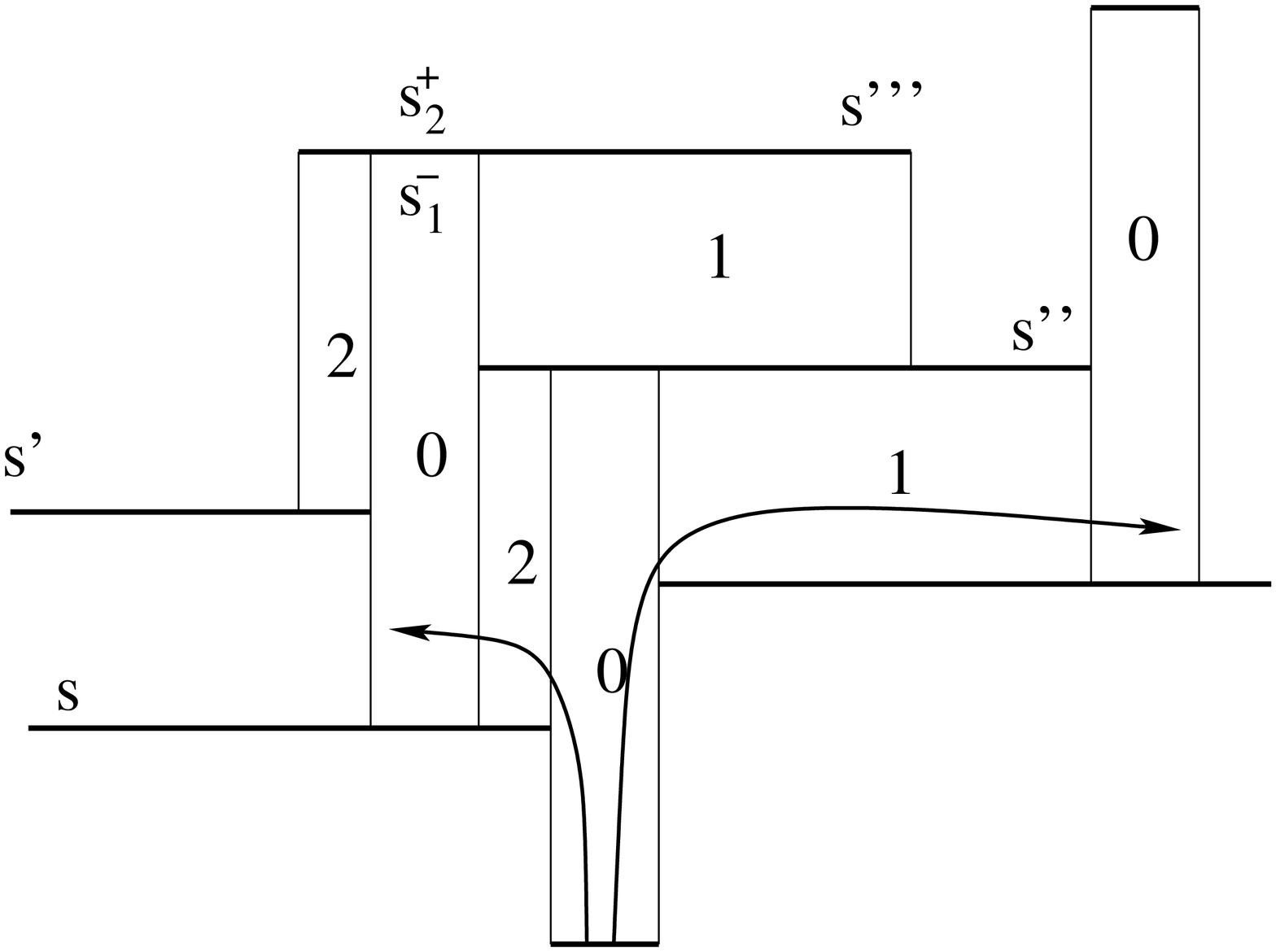}}

$a^m(+):0\rightarrow 2\rightarrow 0$, $b^m(+):0\rightarrow
1\rightarrow 0$.

Fig 13.
\end{center}

\begin{center}
\mbox{\epsfxsize=9cm \epsffile{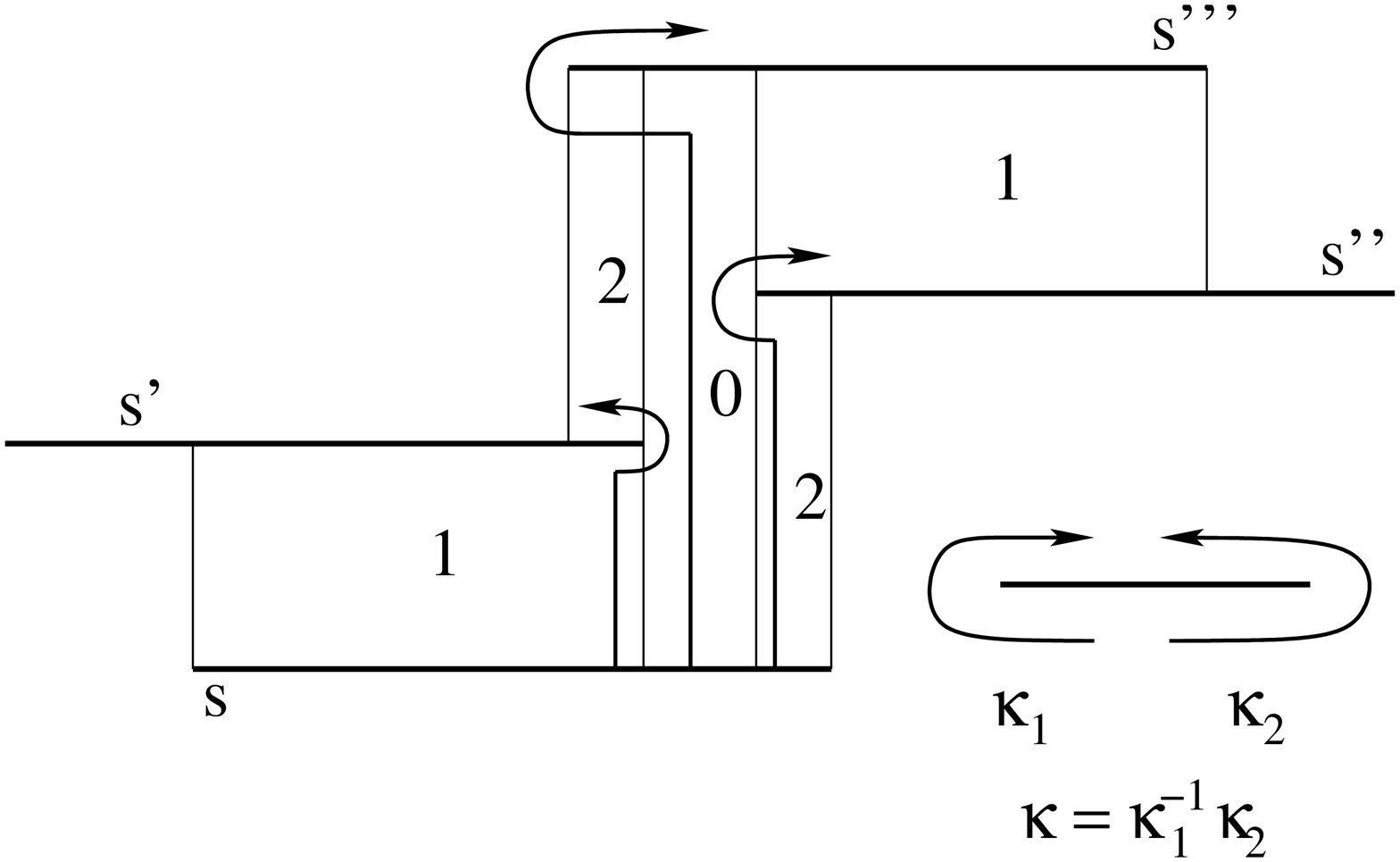}}

Fig 14.
\end{center}

\pagebreak

\begin{center}
\mbox{\epsfxsize=5cm \epsffile{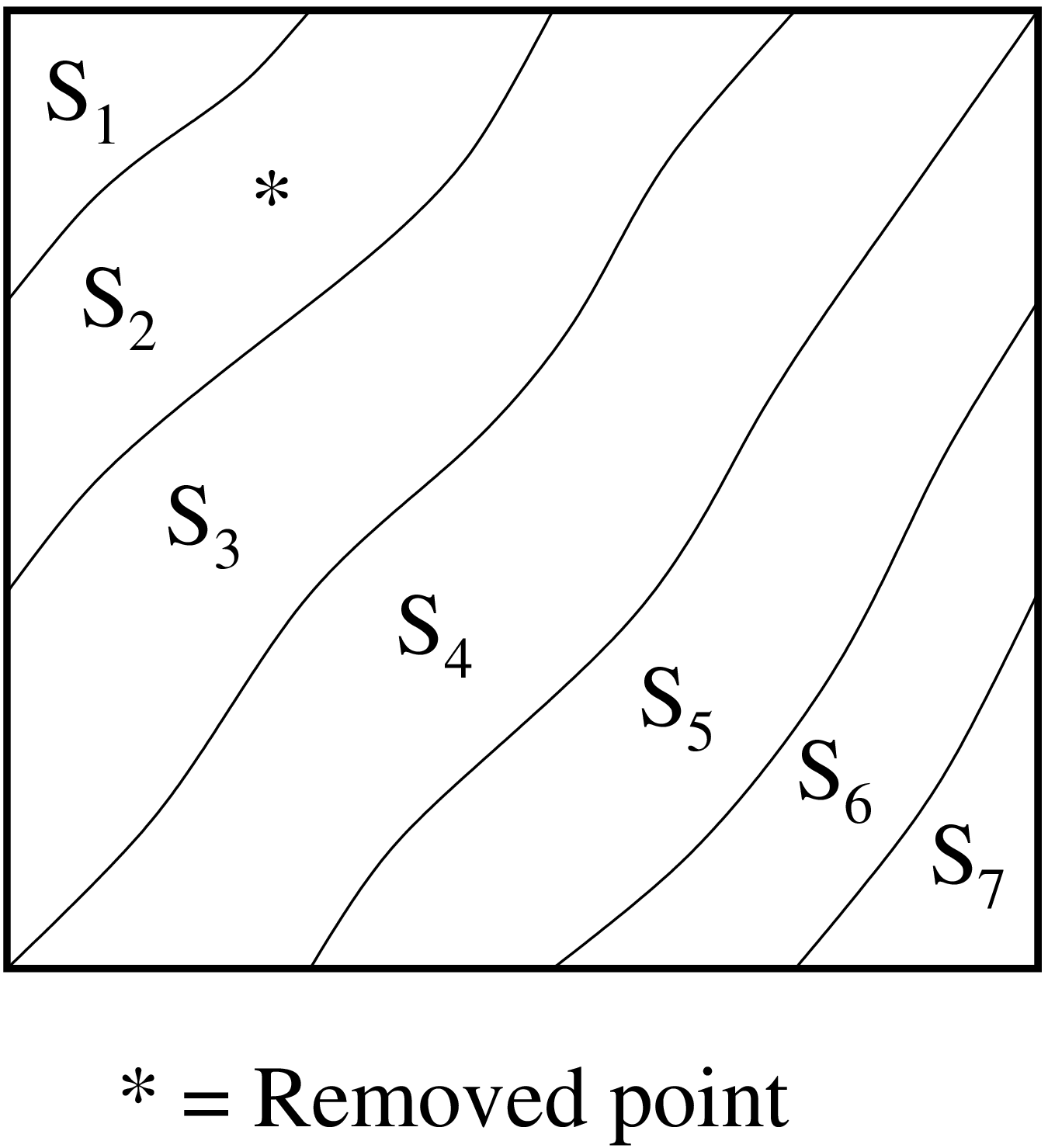}}

$k=4$, $l=3$.

Fig 15.
\end{center}

\begin{center}
\mbox{\epsfxsize=10cm \epsffile{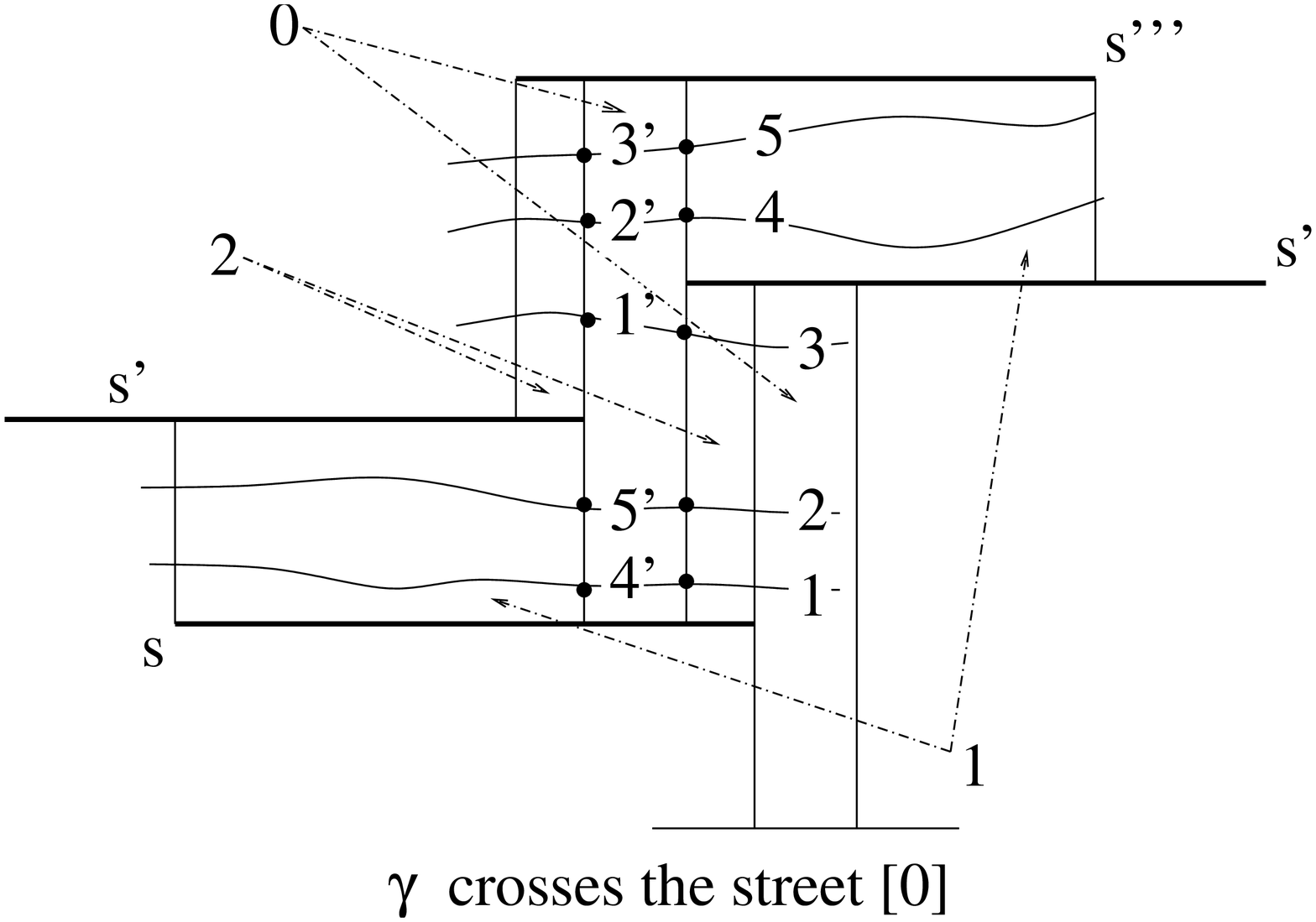}}

$k=3$, $l=2$.

Fig 16.
\end{center}


\vspace{0.3cm}


\begin{thebibliography}{99}

\bibitem{N82}{S.P.Novikov. The Hamiltonian Formalism and Many-Valued Analogs
 of the Morse Theory, Russian Math Surveys
(1982) v 37 n 5 pp 1-56 (in English)}

\bibitem{Abr}{A.A.Abrikosov. Fundamentals of Theory of Metals, North Holland, Amsterdam, 1988
 (in English)}

\bibitem{D99}{I.K.Dynnikov. Geometry of stability regions in the Novikov problem on
 the semiclssical motion of electron, Uspekhi Math Nauk=Russian Math Surveys (1999),
v 54 n 2 pp 21-60 }

\bibitem{MN04}{A.Maltsev, S.Novikov. Dynamical Systems, Topology and Conductivity in Normal Metals.
 Journal of the Statistical Physics, 115, n 1 April 2004,  pp 31-46 }

\bibitem{DN05}{I.Dynnikov, S.Novikov. Topology of  Quasiperiodic Functions on the Plane,
 Uspekhi Math Nauk=Russian Math Surveys (2005) v 60 n 1 }

\bibitem{RdL}{R.Deleo. Numerical Analysis of the Novikov Problem of the motion of electron,
 SIAM Journal For the Applied Dinamical Systems, 2:4, pp 517-545
(2003)}


\bibitem{N65}{S.P.Novikov. Topology of Foliations, Trudy MMO=Transactions of the Moscow
Math Society, 1965, v 14 (in Russian)}

\bibitem{K}{A.~Katok. Invariant measures of flows on oriented surfaces. Soviet
Math. Dokl., v 14, pp 1104--1108  (1973)}

\bibitem{Sat}{
E.~Sataev. The number of invariant measures for flows on
orientable surfaces. Izv. Akad. Nauk SSSR Ser. Mat. v 39, no. 4,
pp 860--878 (1975)}

\bibitem{KHa}{
 B.~Hasselblatt, A.~Katok (ed). Handbook of Dynamical Systems,
Vol. 1A, pp 1015--1089. Elsevier Science B.V. (2002)}



\bibitem{HM}{
J.~Hubbard, H.~Masur. Quadratic differentials and measured
foliations. Acta Math., v 142, pp 221--274 (1979)}

\bibitem{KMS}{
S.~Kerckhoff, H.~Masur,  J.~Smillie. Ergodicity of billiard flows
and quadratic differentials.  Annals of Math., v 124, pp 293--311
(1986)}

\bibitem{M}{
H.~Masur. Interval  exchange transformations and measured
foliations. Ann. of Math., v 115, pp 169-200  (1982)}

\bibitem{V}{
W.~A.~Veech. Gauss measures  for  transformations on the space of
interval exchange maps. Annals  of  Math., v 115, pp 201--242
(1982)}


\bibitem{MT}{
H.~Masur, S.~Tabachnikov. Rational Billiards and Flat Structures.
In: B.~Hasselblatt and A.~Katok (ed), Handbook of Dynamical
Systems, Vol. 1A, 1015--1089. Elsevier Science B.V. (2002)}



\bibitem{Z}{A.Zorich. Proceedings of the Conference ''Geometric Study of
Foliations-Tokyo, 1993/ edited by Muzutani et al, Singapore, World
Scientific (1994)}

\bibitem{Z1}{
A.~Zorich. How do the leaves of a closed 1-form  wind around a
surface.  In the collection: ``Pseudoperiodic Topology'', {\em AMS
Translations}, Ser. 2,  vol.  197, AMS, Providence, RI, pp
135--178 (1999)}

\bibitem{Z2}{A.~Zorich. Flat Surfaces. ``Frontiers in Number Theory, Physics
and Geometry'', Proceedings of Les Houches winter school-2003,
Springer Verlag (2005) (to appear)}

\bibitem{ZK}{
M.~Kontsevich, A.~Zorich: Connected components of the moduli
spaces  of  Abelian differentials.  Invent. Math., 153:3, pp
631--678 (2003)}


\bibitem{Kn}{
M.~Kontsevich. Lyapunov   exponents   and   Hodge  theory. ``The
mathematical beauty of  physics'' (Saclay, 1996), (in Honor of C.
Itzykson) pp 318--332,  Adv.  Ser.  Math.  Phys.,  24. World Sci.
Publishing, River Edge, NJ (1997)}

\bibitem{F}{
G.~Forni. Deviation  of  ergodic  averages  for area-preserving
flows  on surfaces of higher genus.  Annals of Math., v 155, no.
1, pp 1--103  (2002)}




\bibitem{McM}{
C.~McMullen. Dynamics  of  $SL_2(R)$ over moduli space in genus
two. Preprint (2003)}

\bibitem{A}{
V.~I.~Arnold. Topological and ergodic properties of closed 1-forms
with rationally independent periods. Functional Anal. Appl., v 25,
no. 2, pp 81--90 (1991)}

\bibitem{KhS}{
K.~M.~Khanin, Ya.~G.~Sinai. Mixing of some classes of special
flows over rotations of the circle. Functional Anal. Appl. v 26,
no. 3, pp 155--169 (1992)}


\bibitem{L}{G.Levitt, Topology, vol 21 (1982), pp 9-33}

\end{thebibliography}
\end{document}